\documentclass[pdflatex]{amsart}

\title{Higher Nash blowups of the $ A_3 $-singularity}
\author{Rin Toh-yama}

\usepackage{amsmath,amsthm,amssymb}
\usepackage[dvipdfmx]{graphicx}

\usepackage[doublespacing]{setspace}

\theoremstyle{definition}
\newtheorem{Th}{Theorem}[section]
\newtheorem{Prop}[Th]{Proposition}
\newtheorem{Lem}[Th]{Lemma}
\newtheorem{Cor}[Th]{Corollary}
\newtheorem*{MainTh}{Main Theorem}
\newtheorem*{Ques}{Question}
\newtheorem{Def}[Th]{Definition}
\newtheorem{DefProp}[Th]{Definition and Proposition}
\newtheorem{Rem}[Th]{Remark}
\newtheorem{Note}[Th]{Notation}

\newcommand{\GF}[1]{\text{GF}\left(#1\right)}
\newcommand{\interior}[1]{\text{int}\left(#1\right)}
\newcommand{\lm}[1]{\text{lm}_{\preceq}\left( #1 \right)}
\newcommand{\supp}[1]{\text{supp}\left(#1\right)}
\newcommand{\initial}{\text{in}}

\begin{document}
	
	\maketitle
	
	\begin{abstract}
		We show that the $ n $-th Nash blowup of the toric surface singularity of type $ A_3 $ is singular for any $ n > 0 $. It was known that the normalization of the $ n $-th Nash blowup of a toric variety  is also a toric variety associated to the Gr\"{o}bner fan of a certain ideal $ J_n $. In our case, we prove that the Gr\"{o}bner fan contains a non-regular cone. We determine minimal generators of the initial ideal of $ J_n $ with respect to a certain monomial ordering, and show that the reduced Gr\"{o}bner basis of $ J_n $ has polynomials of certain forms for each $ n $.
	\end{abstract}

	\section*{Introduction}
	
	Let  $ X $ be a quasi-projective variety over $ \mathbb{C} $. Classically the Nash blowup of $ X $ was defined in \cite{Nobile}, and recently generalized in \cite{OZ}\cite{Yasuda1} independently.
	
	The classical Nash blowup is defined as follows.
	
	\begin{Def}[\cite{Nobile}]
		Let $ X $ be a subvariety of  $ \mathbb{A}^m $ of $ \dim X = r $ and $ X_{sm} := X \setminus \text{Sing}\left( X \right) $. Let $ G^m_r $ be the Grassmanian of $ r $-dimensional subspaces of $ m $-dimensional $ \mathbb{C} $-space. Then we have the morphism
		\begin{displaymath}
		X_{sm} \hookrightarrow X \times G^m_r;\ P \mapsto \left( P, T_PX \right)
		\end{displaymath}
		where $ T_PX $ is the tangent space of $ X $ at $ P $. Via the morphism, we obtain $ \text{Nash}(X) $  as the closure of $ X_{sm} $ in $ X \times G^m_r $. We also obtain $ \pi :\ \text{Nash}\left( X \right) \rightarrow X $ by restriction of the $ 1 $-st projection $ X \times G^m_r \rightarrow X $. This $ \left( \text{Nash}(X),\pi \right) $ is called the Nash blowup of $ X $. For arbitrary variety, the Nash blowup is defined by gluing of Nash blowups of its affine patches.
	\end{Def}

	This $ \text{Nash}\left( X \right) $ is generalized to the $ n $-th Nash blowup, denoted by $ \text{Nash}_n\left( X \right) $. Let $ \mathcal{M}_{X,P} \subset \mathcal{O}_X $  be an ideal sheaf of a closed point $ P $. In the construction of the classical $ \text{Nash}(X) $, we used the tangent space of each smooth $ P $, which is the dual space of $ \mathcal{M}_{X,P} / \mathcal{M}_{X,P}^2 $. In other words, we considered the $ 1 $-st infinitesimal neighborhood of $ P $, that is the closed subscheme of $ X $ whose ideal sheaf $ \mathcal{M}_{X,P}^2 $. The idea of the generalization is to consider the $ n $-th infinitesimal neighborhood for higher $ n > 0 $.
	
	\begin{Def}[\cite{Yasuda1}]
		Let $ n > 0 $ be an integer.
		\begin{enumerate}
			\item For any closed point $ P \in X $, the $ n $-th fat point $ P^{(n)} $ is defined to be the closed subscheme of $ X $ whose ideal sheaf $ \mathcal{M}^{n+1}_{X,P} $.
			\item Let $ \text{Hilb}_N\left( X \right) $ be the Hilbert scheme of $ N $-points for $ N := \binom{\dim X + n}{\dim X} $. Then $ P \in X_{sm} $ corresponds to $ \left[ P^{(n)} \right] \in \text{Hilb}_N\left( X \right) $, and we have the morphism
			\begin{displaymath}
			X_{sm} \hookrightarrow X \times \text{Hilb}_N\left( X \right);\ P \mapsto \left( P, \left[ P^{(n)} \right] \right).
			\end{displaymath}
			Via the morphism, we obtain $ \text{Nash}_n(X) $  as the closure of $ X_{sm} $ in $ X \times \text{Hilb}_N\left( X \right) $. We also obtain $ \pi_n:\ \text{Nash}_n(X) \rightarrow X $ by restriction of the $ 1 $-st projection $ X \times \text{Hilb}_N\left( X \right) \rightarrow X $. This $ \left( \text{Nash}_n(X),\pi_n \right) $ is called the $ n $-th Nash blowup of $ X $. $ \text{Nash}_1(X) $ was shown to be isomorphic to classical $ \text{Nash}(X) $.
		\end{enumerate}
	\end{Def}
	
	In \cite{Yasuda1} the following question was proposed.
	
	\begin{Ques}[\cite{Yasuda1}, Conjecture 0.2]
		Let $ X $ be any variety of dimension $ d $ and let $ J^{(d-1)} $ be the $ (d-1) $-th neighborhood of the Jacobian subscheme $ J \subset X $ (that is, the closed subscheme defined by $ \mathfrak{j}^d_X $ where $ \mathfrak{j}_X $ is the Jacobian ideal sheaf of $ X $). Let $ [Z] \in \text{Nash}_n(X) $ with $ Z \nsubseteq J^{(d-1)}  $. Then, is $ \text{Nash}_n(X) $ smooth at $ [Z] $?
	\end{Ques}
	
	When $ X $ is a curve, it was proved in \cite{Yasuda1} that $ \text{Nash}_n(X) $ is smooth.
	 
	If the answer were positive for any $ X $, then a resolution of singularities of $ X $ could be obtained as $ \text{Nash}_n(X) $ for $ n \gg 0 $. Hence we could resolve singularities without iterations of operations as with Hironaka's resolution \cite{Hironaka}.
	 
	Our main result shows that such resolutions can not be realized necessarily.
	
	\begin{MainTh}
		Let $ X := \left( z^4 - xy = 0 \right) \subset \mathbb{A}^3 $ be the toric surface singularity of type $ A_3 $. Then $ \text{Nash}_n(X) $ is singular for any $ n > 0 $.
	\end{MainTh}

	The suggestion that the $ A_3 $-singularity might be a counter example was given by T.\ Yasuda (\cite{Yasuda2}, Question 1.4). Moreover extensive calculations were given in \cite{Duarte2} which support the suggestion. We were motivated by these works.
	
	We prove our main theorem as follows. Duarte's theorem (Theorem \ref{GF-theorem}) allows us to describe the normalization $ \overline{\text{Nash}_n(X)} $ of $ \text{Nash}_n(X) $ by using the Gr\"{o}bner fan $ \GF{J_n} $ of a certain ideal $ J_n $. Thus it is enough to see that $ \GF{J_n} $ contains a non-regular cone. On the other hand, we see that maximal cones of $ \GF{J_n} $ come from reduced Gr\"{o}bner bases of $ J_n $. Thus we give the reduced Gr\"{o}bner basis of $ J_n $ with respect to a certain ordering, and show the non-regularity of the cone coming from the basis. Therefore $ \overline{\text{Nash}_n(X)} $ is singular and so is $ \text{Nash}_n(X) $.
	
	This article is organized as follows.
	
	In section $ 1 $, we introduce a general theory of Gr\"{o}bner fans of ideals of monomial subalgebras. The theory is necessary in order to use Duarte's theorem for our result, and we give an explicit way to obtain maximal cones of a Gr\"{o}bner fan from reduced Gr\"{o}bner bases.
	
	In section $ 2 $, we give the proof of our main theorem. First we give a certain monomial ordering $ \preceq $ and determine the minimal generators of $ \initial_{\preceq}\left( J_n \right) $. The determination is the hardest part of this article and needs some ring-theoretic arguments. Then it is shown that the reduced Gr\"{o}bner basis of $ J_n $ with respect to $ \preceq $ has polynomials of certain forms. This enables us to describe the rays of the cone coming from the basis, and the non-regularity is concluded.

	\section{Gr\"{o}bner fans of ideals of monomial subalgebras}
	
	In this section, we introduce a theory of Gr\"{o}bner fans of ideals of monomial subalgebras. A usual theory on Gr\"{o}bner fans is considered for ideals of polynomial rings, but it was known that very analogous results hold for ideals of monomial subalgebras \cite{Duarte1}\cite{Duarte2}\cite{Sturmfels}.
	
	Before we begin to introduce the theory, let us explain why we need it. First let us make clear our settings in this section.
	
	\begin{Note}
		Let $ X $ be any affine toric variety, and $ \sigma \subset \mathbb{R}^d $ be the strongly convex full-dimensional rational polyhedral cone to which $ X $ is associated.
		\begin{enumerate}
			\item $ S := \mathbb{C}[\sigma^\vee \cap \mathbb{Z}^d] $. Then $ X = \text{Spec}\, S $.
			\item Let $ a_1,\ldots,a_s $ generate $ \sigma^\vee \cap \mathbb{Z}^d $, that means $ \sigma^\vee \cap \mathbb{Z}^d = \mathbb{Z}_{\geq 0}a_1 + \cdots + \mathbb{Z}_{\geq 0}a_s $.
			\item  By a coordinate transformation, we can assume $ \sigma^\vee \subset \mathbb{R}_{\geq 0}^d  $. Then $ S $ becomes a monomial subalgebra of $ \mathbb{C}\left[ x_1,\ldots,x_d \right] $ in the following way. For each $ a_i = \left( a_{i,1}, \ldots, a_{i,d} \right) $, take the monomial
			 \begin{displaymath}
			 x^{a_i} := x_1^{a_{i,1}} \cdots x_d^{a_{i,d}} \in \mathbb{C}[x_1,\ldots,x_d].
			 \end{displaymath}
			 Then $ S = \mathbb{C}[x^{a_1},\ldots,x^{a_s}] \subset \mathbb{C}[x_1,\ldots,x_d] $.
		\end{enumerate}
	\end{Note}
	
	With these notations, Duarte's theorem is described as follows.
	
	\begin{Th}[\cite{Duarte1}, Theorem 2.10]\label{GF-theorem}
		$ \overline{\text{Nash}_n(X)} $ is a toric variety. Moreover, let $ J_n := \left\langle \, x^{a_1}-1,\ldots,x^{a_s}-1 \, \right\rangle^{n+1} \subset S $. Then the Gr\"{o}bner fan $ \GF{J_n} $ of $ J_n $ is the fan to which $ \overline{\text{Nash}_n(X)} $ is associated.
	\end{Th}

	By this theorem, we can conclude that $ \overline{\text{Nash}_n(X)} $ is singular if $ \GF{J_n} $ contains a non-regular cone. Now remark that $ J_n $ is an ideal of the monomial subalgebra $ S $, and this is why we need a theory of Gr\"{o}bner fans of ideals of monomial subalgebras.
	
	Now let us begin to introduce the theory. Let $ I $ be arbitrary non-zero ideal of $ S $ till the end of this section.

	\begin{Def}[\cite{Duarte1}, Proposition 1.5]
		Let $ w\in\sigma $.
		\begin{enumerate}		
			\item Let $ 0 \neq f = \sum c_\beta x^\beta \in S $. Put $ m = \max \left\lbrace \, w \cdot \beta \mid x^\beta \in \supp{f} \, \right\rbrace $ where the dot product $ \cdot $ denotes the standard inner product on $ \mathbb{R}^d $. Then we define the initial form of $ f $ with respect to $ w $ to be
			\begin{displaymath}
			\initial_w \left( f \right) := \sum_{
			\substack{ 
				x^\beta \in \supp{f}\\
				w \cdot \beta = m
			}
			} c_{\beta} x^{\beta}.
			\end{displaymath}
			We also define $ \initial_w \left( 0 \right) := 0 $.
			\item $ \initial_w(I):=\langle \, \initial_w(f) \mid f\in I \, \rangle $ is called the initial ideal of $ I $ with respect to $ w $.
			\item Let $ C[w]:=\left\lbrace \, w' \in \sigma \mid \initial_{w'} \left( I \right) = \initial_{w} \left( I \right) \, \right\rbrace  $.
		\end{enumerate}
	\end{Def}

	\begin{DefProp}[\cite{Duarte1}, Definition 1.7]
		Let $ \overline{C[w]} $ be the closure of $ C[w] $ in $ \mathbb{R}^d $. Then
		\begin{displaymath}
		\GF{I}:=\left\lbrace \, \overline{C[w]} \mid w \in \sigma \, \right\rbrace
		\end{displaymath}
		forms a polyhedral fan with $ |\GF{I}| = \sigma $. This is called the Gr\"{o}bner fan of $ I $.
	\end{DefProp}

	We give an alternative description of maximal cones of $ \GF{I} $ below. This is more suitable for our purpose.

	\begin{Def}\label{DefOfCG}
		Let $ \preceq $ be a total ordering on monomials of $ S $. Then $ \preceq $ is a monomial ordering if $ \preceq $ satisfies the following conditions;
		\begin{enumerate}
			\item let $ x^\alpha, x^\beta \in S $. If $ x^\beta $ divides $ x^\alpha $ in $ S $, then $ x^\beta \preceq x^\alpha $;
			\item for any $ x^\gamma \in S $, $ x^\beta \preceq x^\alpha $ implies $ x^{\beta + \gamma} \preceq x^{\alpha + \gamma} $.
		\end{enumerate}
	\end{Def}
	
	\begin{Rem}
		The divisibility between monomials of $ S $ is always considered in not $ \mathbb{C}[x_1,\ldots,x_d] $ but $ S $.
	\end{Rem}
	
	\begin{Def}[\cite{Duarte1}, Definition 1.2, 1.3]
		Let $ \preceq $ be a monomial ordering on $ S $. 
		\begin{enumerate}
			\item A set $ \left\lbrace \, g_1,\ldots,g_t \, \right\rbrace $ of non-zero polynomials of $ I $  is called a Gr\"{o}bner basis of $ I $ with respect to $ \preceq $ if for each $ f \in I \setminus \left\lbrace \, 0 \, \right\rbrace $ there exists $ g_i $ such that $ \lm{g_i} $ divides $ \lm{f} $.
			\item A Gr\"{o}bner basis $ \left\lbrace \, g_1,\ldots,g_t \, \right\rbrace $ is called reduced if, for any $ i $, $ \text{lc}_\preceq\left( g_i \right) = 1 $ and no non-zero monomial of $ g_i $ is divisible by $ \lm{g_j} $ for any $ j \neq i $.
		\end{enumerate} 
	\end{Def}
	
	\begin{Th}[\cite{Duarte1}, Theorem 1.4]\label{ExistenceOfReducedGB}
		Let $ \preceq $ be a monomial ordering on $ S $. Then $ I $ has a unique reduced Gr\"{o}bner basis with respect to $ \preceq $.
	\end{Th}

	\begin{Def}\label{CG-definition}
			\begin{enumerate}
				\item Let $ \left\lbrace \, g_1,\ldots,g_t \, \right\rbrace $ be the reduced Gr\"{o}bner basis of $ I $ with respect to $ \preceq $. Then
				\begin{displaymath}
				\mathbb{G} := \left\lbrace \, \left( g_1,\lm{g_1} \right),\ldots,\left( g_t,\lm{g_t} \right) \, \right\rbrace
				\end{displaymath}
				is called the marked Gr\"{o}bner basis of $ I $ with respect to $ \preceq $. The marked Gr\"{o}bner basis for a monomial ordering is referred to as `` a marked Gr\"{o}bner basis ''.
				\item Let $ \mathbb{G} = \left\lbrace \, \left( g_1,x^{\alpha_1} \right),\ldots,\left( g_t,x^{\alpha_t} \right) \, \right\rbrace $ be a marked Gr\"{o}bner basis of $ I $. Then the cone $ C_\mathbb{G} \subset \sigma $ is defined to be
				\begin{displaymath}
				C_\mathbb{G} := \left\lbrace \, w\in\sigma \mid (\alpha_i-\beta)\cdot w\geq 0 \text{ for any }\alpha_i \text{ and } x^\beta \in \supp{g_i} \, \right\rbrace.
				\end{displaymath}
			\end{enumerate}
	\end{Def}

	We see that maximal cones of $ \GF{I} $ are exactly cones given as $ C_\mathbb{G} $.
	
	\begin{Lem}\label{Robbiano'sTheorem}
		Let $ \preceq $ be a monomial ordering on $ S $, which we regard as a total ordering on $ \sigma^\vee \cap \mathbb{Z}^d $.
		
		Then there exists vectors $ w_1,\ldots,w_r \in \mathbb{R}^d $ such that $ \preceq $ extends to a total ordering on $ \mathbb{Q}^d \subset \mathbb{R}^d $ as follows. For any $ \alpha, \beta \in \mathbb{Q}^d $, $ \beta \preceq \alpha $ if and only if there exists $ r_0 \leq r $ such that
		\begin{displaymath}
		\forall i < r_0,\ \left( \alpha - \beta \right) \cdot w_i = 0\ \text{and}\ \left( \alpha - \beta \right) \cdot w_{r_0} > 0.
		\end{displaymath}
		In this case, we say that $ \preceq $ is associated to the $ r \times d $ matrix whose $ i $-th row is $ w_i $.
		
		\begin{proof}
			Let $ H \subset \mathbb{Z}^d $ be the abelian subgroup generated by $ \sigma^\vee \cap \mathbb{Z}^d $. Then $ \mathbb{Q}\otimes_\mathbb{Z} H = \mathbb{Q}^d $ since $ \sigma^\vee \subset \mathbb{R}^d $ is full-dimensional.
			
			One can easily check that $ \preceq $ extends to a total ordering on $ H $ as follows. For any $ p, p' \in H $, take expressions $ p = p_+ - p_- $ and $ p' = p'_+ - p'_- $ for some $ p_+, p_-, p'_+, p'_- \in \sigma^\vee \cap \mathbb{Z}^d $. Then $ p \preceq p' $ if and only if $ p_+ + p'_- \preceq p'_+ + p_- $.
			
			Moreover $ \preceq $ extends to a total ordering on $ \mathbb{Q}^d = \mathbb{Q}\otimes_\mathbb{Z} H $ as follows. For any $ q, q' \in \mathbb{Q}\otimes_\mathbb{Z} H $, there exists $ r \in \mathbb{Z}_{> 0} $ such that $ rq,rq'\in H $. Then $ q \preceq q' $ if and only if $ rq \preceq rq' $.
			
			Now Robbiano's theorem (\cite{Rob}, Theorem 4) shows that there exists $ r > 0 $ and a real $ r \times d $ matrix $ M $ such that the ordering $ \preceq $ on $ \mathbb{Q}^d $ is associated to $ M $. Then $ w_i := (\, i\text{-th row of}\ M \, ) $ are the expected ones.
		\end{proof}
	\end{Lem}
	
	\begin{Lem}[c.\ f.\ \cite{CLO}, Chapter 8, Theorem 4.7]\label{CG-lemma}
		Let $ \mathbb{G} = \left\lbrace \, \left( g_1,x^{\alpha_1} \right),\ldots,\left( g_t,x^{\alpha_t} \right) \, \right\rbrace $ be a marked Gr\"{o}bner basis of $ I $. Then
		\begin{enumerate}
			\item $ C_\mathbb{G} $ is a strongly convex full-dimensional rational polyhedral cone.
			
			\item For any $ w \in \sigma $, $ w \in \interior{C_\mathbb{G}} $ if and only if
			\begin{displaymath}
			(\alpha_i-\beta)\cdot w > 0 \text{ for any } i \text{ and } x^\beta \in \supp{g_i}\setminus \left\lbrace \, x^{\alpha_i} \, \right\rbrace.
			\end{displaymath}
		\end{enumerate}
		
		\begin{proof}
			(1) $ C_\mathbb{G} $ is a rational polyhedral cone because the entries of $ \alpha_i-\beta $ are rational. Moreover $ C_\mathbb{G} $ is strongly convex since $ \sigma \supset C_\mathbb{G} $ is strongly convex.
			
			Let us show the full-dimensionality of $ C_\mathbb{G} $. By the definition of $ C_\mathbb{G} $, it is clear that $ C_\mathbb{G} $ contains an open subset $ U $ of $ \mathbb{R}^d $ defined as
			\begin{displaymath}
			U:= \interior{\sigma} \cap \bigcap_{
				1 \leq i \leq t} \left\lbrace \, w \in \mathbb{R}^d \mid \left(\alpha_i - \beta \right) \cdot w > 0 \ \text{for all}\ x^\beta \in \supp{g_i}\setminus \left\lbrace \, x^{\alpha_i} \, \right\rbrace\, \right\rbrace
			\end{displaymath}
			Therefore it is enough to see $ U \neq \emptyset $.
			
			Let $ \preceq $ be a monomial ordering on $ S $ which provides $ \mathbb{G} $. By Lemma \ref{Robbiano'sTheorem}, $ \preceq $ extends to a total ordering on $ \mathbb{Q}^d $ associated to some $ r \times d $ matrix $ M $. Let $ w_i $ be the $ i $-th row of $ M $ and put $ w_\epsilon := w_1+\epsilon w_2+\cdots + \epsilon^{r-1}w_r \in \mathbb{R}^d $. We will see $ w_\epsilon \in U $ for sufficiently small $ \epsilon > 0 $.
			
			First, let us observe the following. Fix any $ \gamma_1, \gamma_2 \in \sigma^\vee \cap \mathbb{Z}^d $ with $ x^{\gamma_1} \prec x^{\gamma_2} $. Then there exists $ r_0 \leq r $ such that $ \left( \gamma_2 - \gamma_1 \right) \cdot w_i = 0\ $ for all $ i <r_0 $ and $\ \left( \gamma_2 - \gamma_1 \right) \cdot w_{r_0} > 0 $. Therefore, for sufficiently small $ \epsilon > 0 $, we have  $ (\gamma_2 - \gamma_1)\cdot \left( w_{r_0} + \epsilon w_{r_0+1}+\cdots + \epsilon^{r-r_0}w_r \right) > 0 $. This implies $ (\gamma_2 - \gamma_1) \cdot w_\epsilon = (\gamma_2 - \gamma_1) \cdot \left( w_1+\epsilon w_2+\cdots + \epsilon^{r-1}w_r \right) > 0 $.
			
			Let $ L $ be any ray of $ \sigma^\vee $ and $ \mu_L $ be a lattice point of $ \text{rel.int}\left(L\right) $. Now $ 1 \prec x^{\mu_L} $. Then, by above observation, we have $ \mu_L \cdot w_\epsilon > 0 $ for sufficiently small $ \epsilon > 0 $. By restricting $ \epsilon $ for all $ L $, we have $ w_\epsilon \in \interior{\sigma} $. Moreover, for any $ x^\beta \in \supp{g_i}\setminus \left\lbrace \, x^{\alpha_i} \, \right\rbrace  $, we have $ x^\beta \prec x^{\alpha_i} $. Thus $ (\alpha_i - \beta) \cdot w_\epsilon > 0 $ for sufficiently small $ \epsilon > 0 $.
			
			Therefore $ w_\epsilon \in U $ for sufficiently small $ \epsilon > 0 $, and hence (1) holds.
			
			(2) follows from the full-dimensionality of $ C_\mathbb{G} $.
		\end{proof}
	\end{Lem}

	\begin{Lem}[\cite{Duarte2}, Appendix A, Proposition A.2.2]\label{initial-lemma} Let $ \preceq, \preceq' $ be monomial orderings on $ S $.
		Then $ \initial_{\preceq'}\left( I \right) \subset \initial_{\preceq}\left( I \right) $ implies $ \initial_{\preceq'}\left( I \right) = \initial_{\preceq}\left( I \right) $.
	\end{Lem}
	
	\begin{Def}
		Let $ \preceq $ be a monomial ordering on $ S $ and $ w \in \sigma $. Then $ w $-weighted ordering associated to $ \preceq $, denoted by $ \preceq_w $, is defined as follows;
		\begin{displaymath}
		x^\beta \preceq_w x^\alpha \Leftrightarrow \left( \left( \alpha - \beta \right) \cdot w > 0 \right) \text{ or } \left( \left( \alpha - \beta \right) \cdot w = 0 \text{ and } x^\beta \preceq x^\alpha \right).
		\end{displaymath}
		One can easily see that $ \preceq_w $ is also a monomial ordering on $ S $.
	\end{Def}
	
	\begin{Lem}[\cite{Duarte2}, Appendix A, proof of Proposition A.3.1]\label{C[w]-lemma}
		Let $ w \in \sigma $ and $ \preceq $ be any monomial ordering on $ S $. Let $ G $ be the reduced Gr\"{o}bner basis of $ I $ with respect to $ \preceq_w $. Then 
		\begin{displaymath}
		C[w]=\left\lbrace \, w'\in\sigma \mid \initial_{w'}(g)=\initial_{w}(g)\ \text{for all}\ g \in G \, \right\rbrace.
		\end{displaymath}
	\end{Lem}

	\begin{Cor}\label{G-C[w]-corollary}
		Let $ \mathbb{G} = \left\lbrace \, \left(g_1, x^{\alpha_1}\right),\ldots,\left( g_t, x^{\alpha_t} \right) \, \right\rbrace $ be a marked Gr\"{o}bner basis of $ I $ and $ w \in C_\mathbb{G} $.
		\begin{enumerate}
			\item Let $ \preceq $ be any monomial ordering on $ S $ which provides $ \mathbb{G} $. Then $ \preceq_w $ also provides $ \mathbb{G} $.
			\item If $ w \in \interior{C_\mathbb{G}} $, then $ C[w] = \interior{C_\mathbb{G}} $.
		\end{enumerate}
		
		\begin{proof}
			(1) It is enough to see
			\begin{displaymath}
			\text{lm}_{\preceq_{w}}\left( g_i \right) = x^{\alpha_i}\ \text{for}\ 1 \leq i \leq t\ \text{and}\ \initial_{\preceq_{w}}\left( I \right) = \langle \,  x^{\alpha_1},\ldots,x^{\alpha_t}\, \rangle.
			\end{displaymath}
			By the definition of $ C_\mathbb{G} $, $ w \in C_\mathbb{G} $ satisfies $ \lm{g_i} = x^{\alpha_i} \in \supp{\initial_{w}\left( g_i \right)} $. Then $ \text{lm}_{\preceq_{w}}\left( g_i \right) = \lm{\initial_{w}\left( g_i \right)} = x^{\alpha_i} $. Now we have
			\begin{displaymath}
			\initial_{\preceq}\left( I \right) = \langle \, x^{\alpha_1},\ldots,x^{\alpha_t}\, \rangle \subset \initial_{\preceq_{w}}\left( I \right).
			\end{displaymath}
			 Hence $ \initial_{\preceq}\left( I \right) = \initial_{\preceq_{w}}\left( I \right) $ by Lemma \ref{initial-lemma}, and hence the assertion holds.
			
			(2) By (1) $ \mathbb{G} $ is the marked Gr\"{o}bner basis of $ I $ with respect to $ \preceq_w $ and Lemma \ref{CG-lemma} (2) shows that $ \initial_{w}\left(g_i\right) = x^{\alpha_i} $ for any $ 1 \leq i \leq t $. Then, by Lemma \ref{C[w]-lemma},
			\begin{displaymath}
			C[w]=\left\lbrace \, w'\in\sigma \mid \initial_{w'}(g_i)=x^{\alpha_i}\ \text{for all}\ 1 \leq i \leq t \, \right\rbrace. 
			\end{displaymath}
			Hence $ C[w] = \interior{C_\mathbb{G}} $ by Lemma \ref{CG-lemma} (2).
		\end{proof}
	\end{Cor}

	\begin{Lem}\label{w1-lemma}
		Let $ \preceq $ be a monomial ordering associated to a matrix $ M $ and $ \mathbb{G} $ be the marked Gr\"{o}bner basis of $ I $ with respect to $ \preceq $. Then $ w_1 \in C_\mathbb{G} $ where $ w_1 $ is the $ 1 $-st row of $ M $.
		\begin{proof}
			The assertion follows from the definitions of $ C_\mathbb{G} $ and $ \preceq $.
		\end{proof}
	\end{Lem}

	\begin{Th}
		If $ \mathbb{G} $ is a marked Gr\"{o}bner basis of $ I $, then $ C_\mathbb{G} $ is an element of $ \GF{I} $. Conversely, any maximal cone of $ \GF{I} $ is given as $ C_\mathbb{G} $ for some $ \mathbb{G} $.
	
		\begin{proof}
			$ C_\mathbb{G} \in \GF{I} $ follows from Corollary \ref{G-C[w]-corollary} (2). Conversely, fix any maximal cone $ C \in \GF{I} $ and take $ w \in \interior{C} $. Let $ \preceq $ be any monomial ordering on $ S $ and $ \mathbb{G} $ be the marked Gr\"{o}bner basis of $ I $ with respect to $ \preceq_w $. Then $ w \in C_{\mathbb{G}} $ by Lemma \ref{w1-lemma}, and $ C_{\mathbb{G}} $ is an element of $ \GF{I} $. This implies that $ C $ is a face of $ C_\mathbb{G} $. By maximality of $ C $, we have $ C = C_\mathbb{G} $.
		\end{proof}
	\end{Th}

	\section{Higher Nash blowups of the $ A_3 $-singularity}
	
	We give the proof of our main theorem.
	
	\begin{Note}\label{Settings} In this section let $ X := \left(z^4-xy = 0\right) \subset \mathbb{A}^3 $.
		\begin{enumerate}
			\item Let $ \sigma \subset \mathbb{R}^2 $ be the cone generated by $ \left(0,1 \right) $ and $ \left( 4,-3 \right) $. The dual cone $ \sigma^\vee \subset \mathbb{R}^2 $ is generated by $ \left(1,0\right) $ and $ \left( 3,4 \right) $. The both cones are strongly convex and full-dimensional.
			\item The semi-group $ \sigma_{\mathbb{Z}} := \sigma^\vee \cap \mathbb{Z}^2 $ is generated by $ \left(1,0\right),\left(3,4\right),\left(1,1\right) $.
			\item $ S := \mathbb{C}[\sigma_{\mathbb{Z}}] = \mathbb{C}[u,u^3v^4,uv] \subset \mathbb{C}[u,v] $. There is a surjective morphism
			\begin{displaymath}
			F:\ \mathbb{C}[x,y,z] \twoheadrightarrow S;\ x \mapsto u,\ y\mapsto u^3v^4,\ z\mapsto uv
			\end{displaymath}
			with $ \ker F = \langle \, z^4-xy \, \rangle $. Hence $ X $ is isomorphic to $ \text{Spec}\, S $, the toric variety associated to $ \sigma $.
			\item For any integer $ n > 0 $, we put
			\begin{displaymath}
			J_n := \langle \, u-1,u^3v^4-1,uv-1 \, \rangle^{n+1} \subset S.
			\end{displaymath}
			Then the normalization $ \overline{\text{Nash}_n(X)} $ of $ \text{Nash}_n(X) $ is the toric variety associated to $ \GF{J_n} $ (Theorem \ref{GF-theorem}).
		\end{enumerate}
	\end{Note}

	\begin{Rem}
		Points of $ \sigma_{\mathbb{Z}} $ correspond to monomials in $ S $ bijectively, hence they are identified without explicit notice. For example, we identify $ \left( 3,4 \right) \in \sigma_{\mathbb{Z}} $ with $ u^3v^4 \in S $, and $ (1,0) + (1,1) $ with $ u\cdot uv $.
	\end{Rem}

	We find a non-regular cone in $ \GF{J_n} $ to conclude that $ \overline{\text{Nash}_n(X)} $ is singular. As we explained in the previous section, it is enough to find a marked Gr\"{o}bner basis $ \mathbb{G}_n $ of $ J_n $ such that $ C_{\mathbb{G}_n} $ is non-regular.
	
	\begin{Def}\label{<-definition}
		\begin{enumerate}
			\item Let $ \preceq $ be the monomial ordering on $ S $ associated to
			\begin{displaymath}
			\left(\begin{array}{cc}
			2 & -1 \\
			1 & 1 
			\end{array}
			\right).
			\end{displaymath}
			\item Let $ \mathbb{G}_n $ be the marked Gr\"{o}bner basis  of $ J_n $ with respect to $ \preceq $.
			\item Let $ \mathbb{M}_n := \left\lbrace \, \alpha \mid \left( g, \alpha \right) \in \mathbb{G}_n \, \right\rbrace $.
		\end{enumerate}
	\end{Def}

	We describe $ C_{\mathbb{G}_n} $ explicitly and see the non-regularity. We first study $ \mathbb{M}_n $.

	\subsection{Candidate for $ \mathbb{M}_n $}

\begin{figure}[htbp]
		\begin{tabular}{cc}
		\begin{minipage}{0.5\hsize}
			\centering
			\includegraphics[keepaspectratio, width=\hsize]{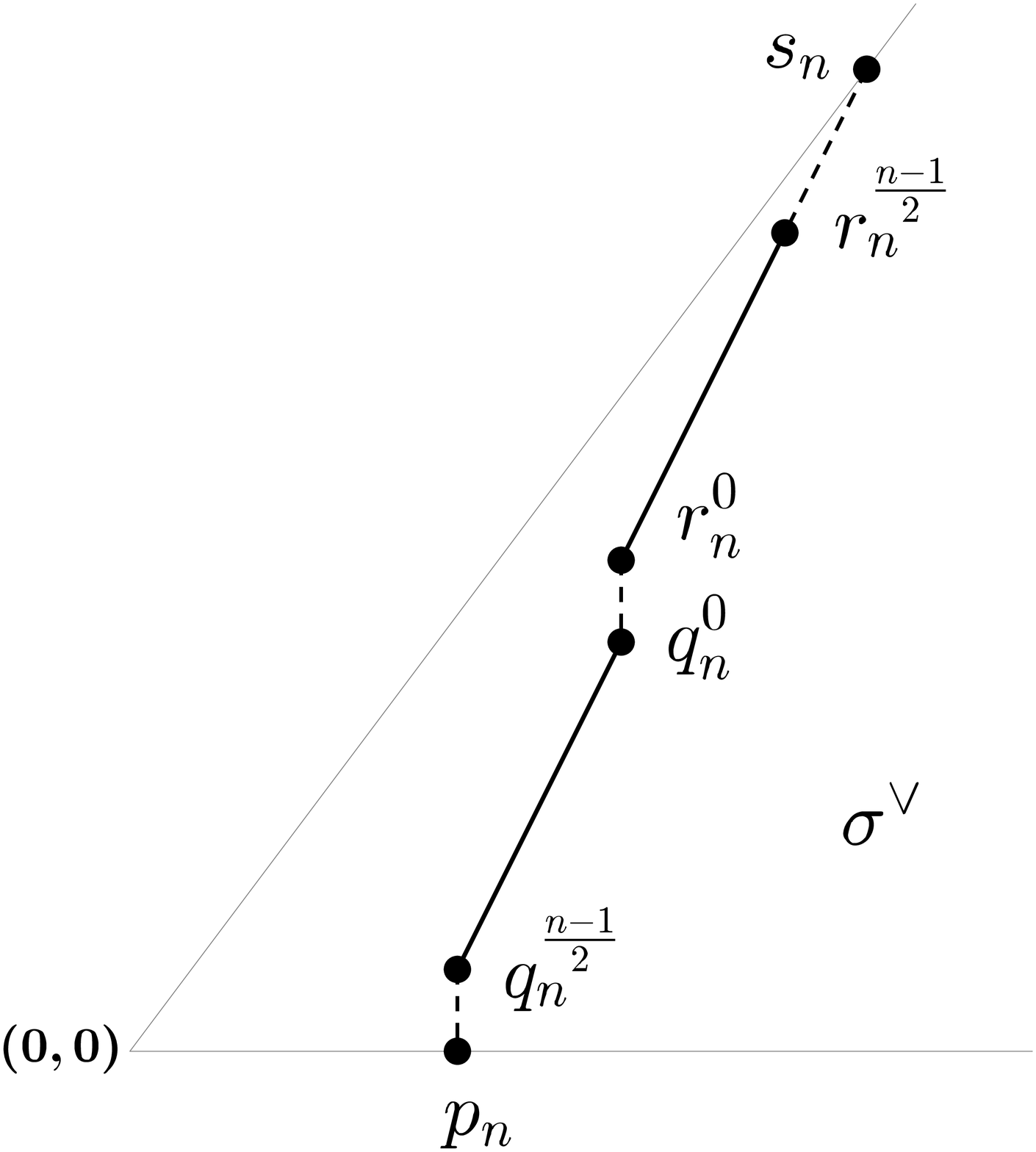}
			\caption{$ \mathcal{P}_n $ for odd $ n $} \label{Pn-odd}
		\end{minipage}
		\begin{minipage}{0.5\hsize}
			\centering
			\includegraphics[keepaspectratio, width=\hsize]{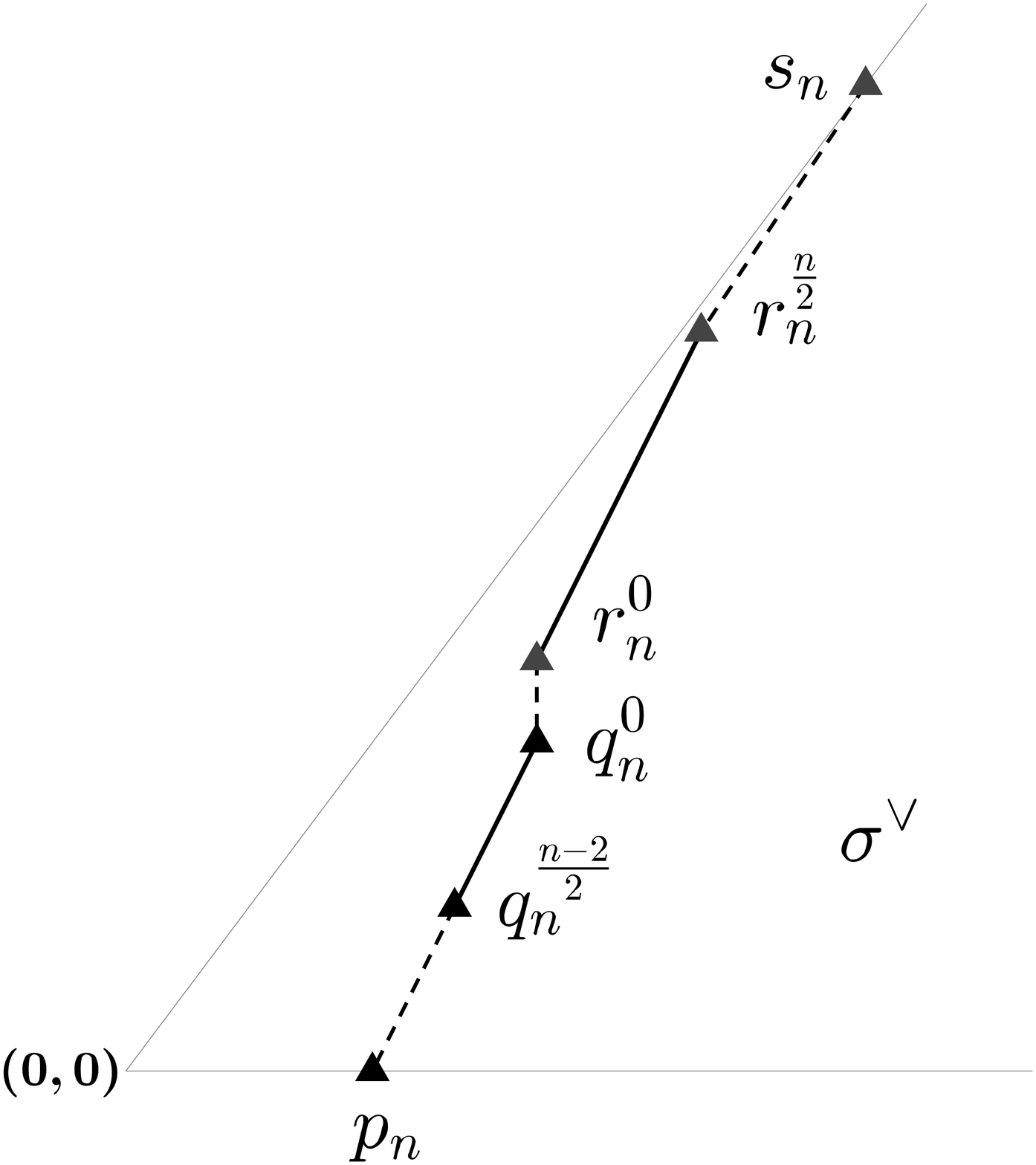}
			\caption{$ \mathcal{P}_n $ for even $ n $}\label{Pn-even}
		\end{minipage}
	\end{tabular}
	\end{figure}

	\begin{Def}\label{Pn-definition}
		For each integer $ n > 0 $, let $ \mathcal{P}_n $ be the set consisting of following points of $ \sigma_\mathbb{Z} $ (see Figure \ref{Pn-odd} and Figure \ref{Pn-even}); if $ n $ is odd, then 
		\begin{displaymath}
		\begin{array}{ll}		
		p_n := \left( \frac{n+3}{2},0 \right),& \\
		q_n^0 := \left( \frac{n+3}{2},1 \right)+\frac{n-1}{2}\left( 1,2 \right), & \ q_n^i := q_n^0 - i\left( 1,2 \right)\ \left( 0 \leq i \leq \frac{n-1}{2} \right),\\
		r_n^0 := q_n^0 + (0,1), & \ r_n^j := r_n^0 + j\left( 1,2 \right)\ \left( 0 \leq j \leq \frac{n-1}{2} \right),\\
		s_n := \frac{n+1}{2}\left( 3,4 \right);&
		\end{array}
		\end{displaymath}
		if $ n $ is even, then
		\begin{displaymath}
		\begin{array}{ll}		
		p_n := \left( \frac{n+2}{2},0 \right), & \\
		q_n^0 := \left( \frac{n+2}{2},0 \right)+\frac{n}{2}\left( 1,2 \right), & q_n^i := q_n^0 - i\left( 1,2 \right)\ \left( 0 \leq i \leq \frac{n-2}{2} \right),\\
		r_n^0 := q_n^0 + (0,1), & r_n^j := r_n^0 + j\left( 1,2 \right)\ \left( 0 \leq j \leq \frac{n}{2} \right),\\
		s_n := \left( \frac{n+2}{2} \right)\left( 3,4 \right). &
		\end{array}
		\end{displaymath}
	\end{Def}

	We see $ \mathbb{M}_n = \mathcal{P}_n $ below.
	
	\begin{Lem}\label{innerPn-lemma}
		Let $ n > 0 $ be an integer.
		\begin{enumerate}
			\item If $ n $ is odd, then $ p_n = q_n^{\frac{n-1}{2}} - (0,1),\ s_n = r_n^{\frac{n-1}{2}} + (1,2) $.
			\item If $ n $ is even, then $ p_n = q_n^\frac{n-2}{2} - (1,2),\ s_n = r_n^{\frac{n}{2}} + (2,3) $.
		\end{enumerate}
		\begin{proof}
			The assertions follow from direct calculations.
		\end{proof}
	\end{Lem}

	Let us describe $ \mathcal{P}_n $. By the definition of $ \mathcal{P}_n $ and Lemma \ref{innerPn-lemma}, we obtain Figure \ref{Pn-odd} and Figure \ref{Pn-even}; all $ q_n^i $ and $ r_n^j $ are lying on the segments of a thick line, and conversely all lattice points on the segments are members of $ \mathcal{P}_n $. The segments of a broken line have lattice points only at the edges.
	
	If $ n $ is odd, then the slopes of segments $ q_n^\frac{n-1}{2}q_n^0 $ and $ r_n^0s_n $ are same $ 2 $. If $ n $ is even, then the slopes of $ p_nq_n^0 $ and $ r_n^0r_n^\frac{n}{2} $ are same $ 2 $, and the slope of $ r_n^\frac{n}{2}s_n $ is $ \frac{3}{2} $.
	
	\begin{Lem}\label{Pn+sigma-lemma}
		In Figure \ref{Pn-odd} and Figure \ref{Pn-even}, there is a polygonal line over $ \mathcal{P}_n $ dividing $ \sigma^\vee $ into two regions. Then $ \mathcal{P}_n + \sigma_\mathbb{Z} $ consists of all lattice points of the region on the right side.
		
		\begin{proof}
			Let $ U $ be the region on the right side. We have
			\begin{displaymath}
			\mathcal{P}_n + \sigma_\mathbb{Z} = \bigcup_{a \in \mathcal{P}_n} \left(a + \sigma_\mathbb{Z} \right).
			\end{displaymath}
			For any $ a \in \mathcal{P}_n $, it is clear that $ a + \sigma_\mathbb{Z} $ consists of lattice points of $ a + \sigma^\vee $. Figure \ref{rn0+sigmaforodd} and Figure \ref{rn0+sigmaforeven} describe $ a + \sigma^\vee $ for $ a = r_n^0 $, thus one can see that $ a + \sigma_\mathbb{Z} $ is contained in $ U $; indeed the ray of $ \sigma^\vee $ generated by $ (3,4) $ has the slope $ \frac{4}{3} < 2, \frac{3}{2} $, and hence $ a $ is the only point of the polygonal line contained in $ a + \sigma^\vee $.
			
			\begin{figure}[htbp]
				\begin{tabular}{cc}
					\begin{minipage}{0.5\hsize}
						\centering
						\includegraphics[keepaspectratio, width=\hsize]{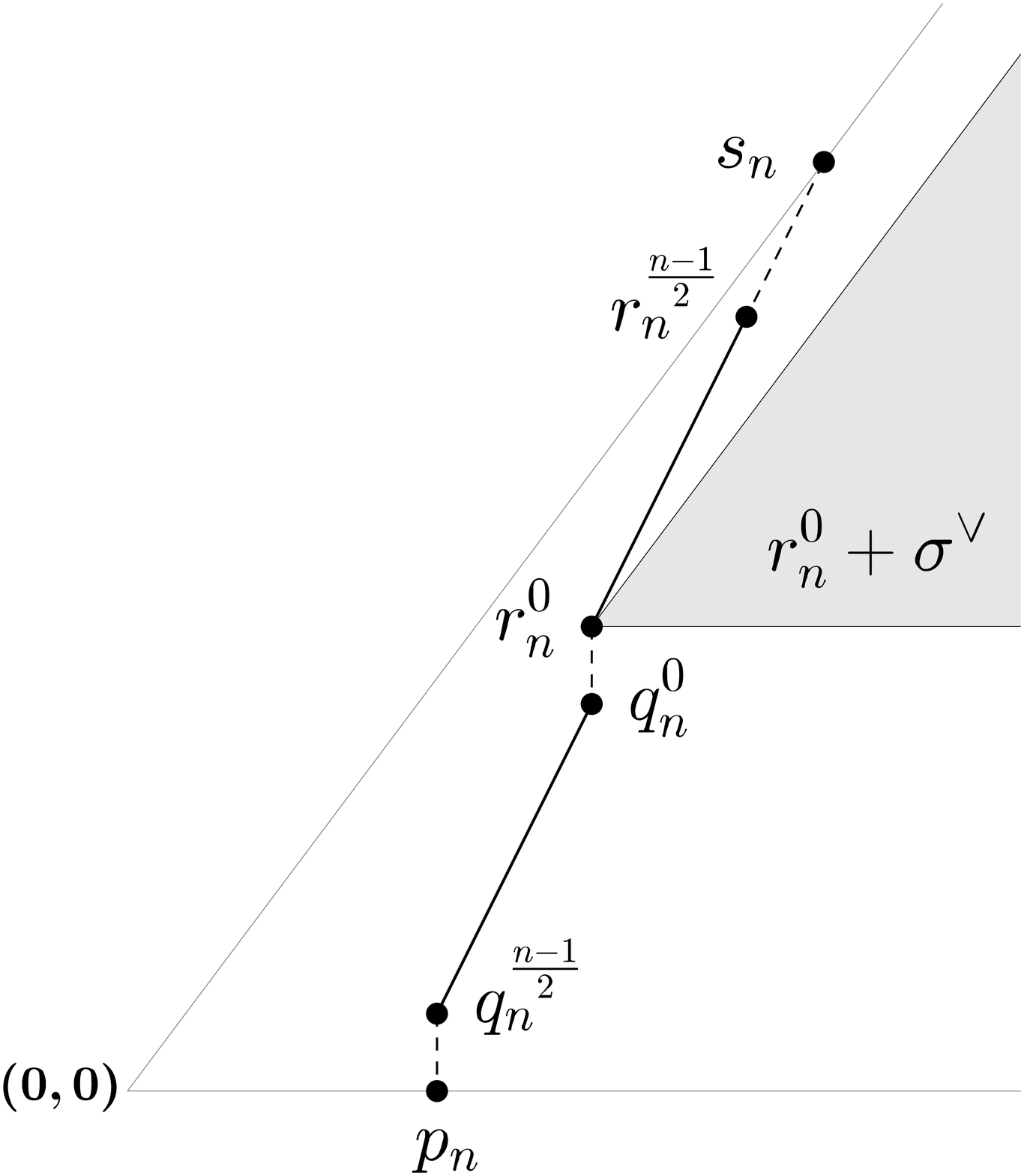}
						\caption{$ r_n^0 + \sigma^\vee $ for odd $ n $} \label{rn0+sigmaforodd}
					\end{minipage}
					\begin{minipage}{0.5\hsize}
						\centering
						\includegraphics[keepaspectratio, width=\hsize]{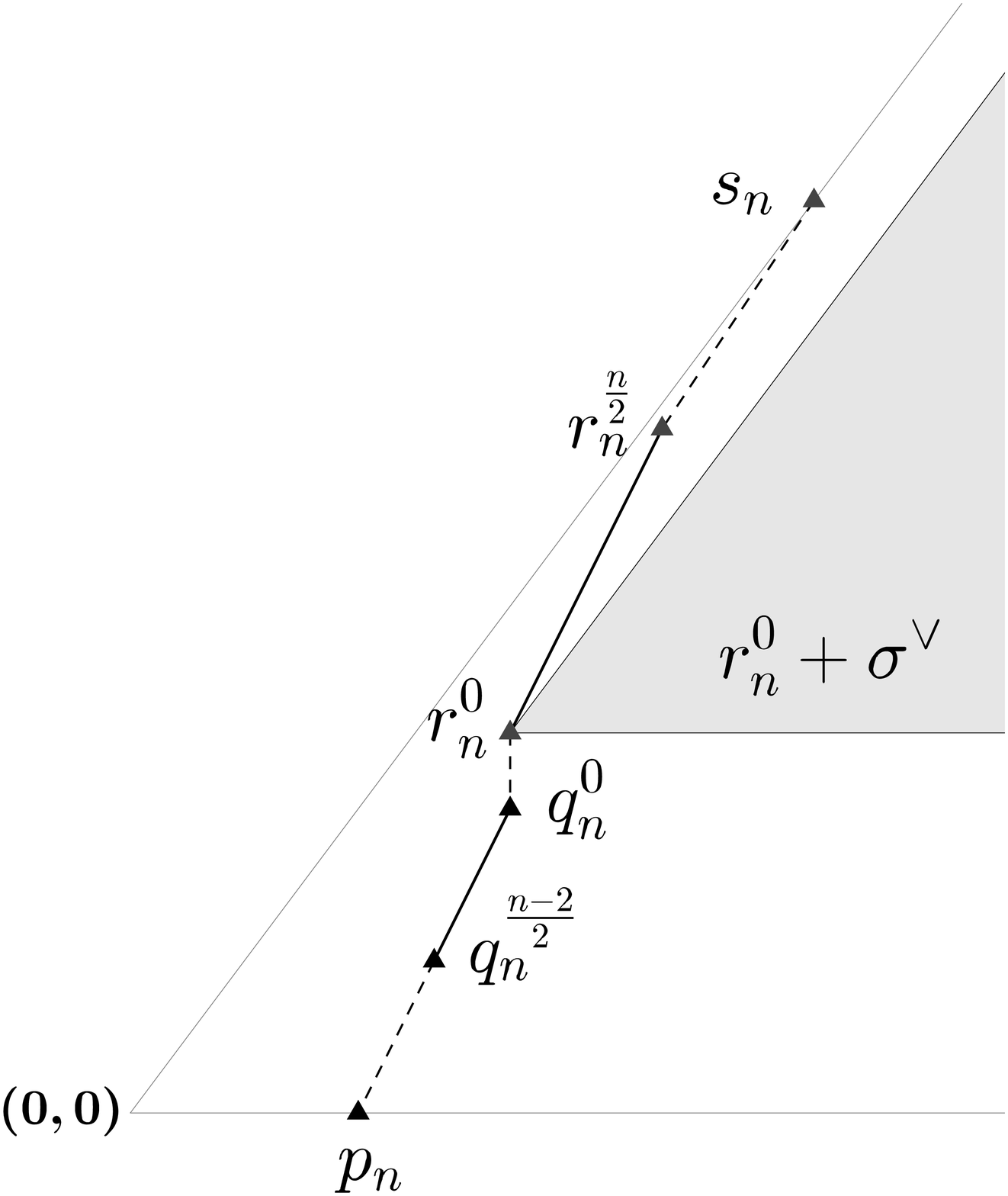}
						\caption{$ r_n^0 + \sigma^\vee $ for even $ n $}\label{rn0+sigmaforeven}
					\end{minipage}
				\end{tabular}
			\end{figure}

			Therefore $ \mathcal{P}_n + \sigma_\mathbb{Z} $ is contained in $ U $.
			
			Take any lattice point $ b \in U $, and let us show $ b \in \mathcal{P}_n + \sigma_\mathbb{Z} $.
			
			Let $ L $ be the line passing through both $ b $ and $ b + (1,0) $, and $ b_0 $ be the lattice point of $ U $ such that $ b_0 \in L $ and $ b_0 - (1,0) \notin U $. Then it is enough to see $ b_0 \in \mathcal{P}_n + \sigma_\mathbb{Z} $; indeed, in this case, we have $ b \in b_0 + \sigma_\mathbb{Z} \subset \mathcal{P}_n + \sigma_\mathbb{Z} $.
			
			Now we can assume $ b \notin s_n + \sigma_\mathbb{Z} $, then $ L $ has an intersection $ c $ with the polygonal line. Figure \ref{nbhofb0} describes this situation when $ c $ lies on the segment $ aa' $ of the slope $ 2 $ for some $ a,a' \in \mathcal{P}_n $. Then one can see $ b_0 \in a + \sigma_\mathbb{Z} $. The other cases are easily checked by similar figures.
			
			\begin{figure}[htbp]
				\centering
				\includegraphics[keepaspectratio, width=0.5\hsize]{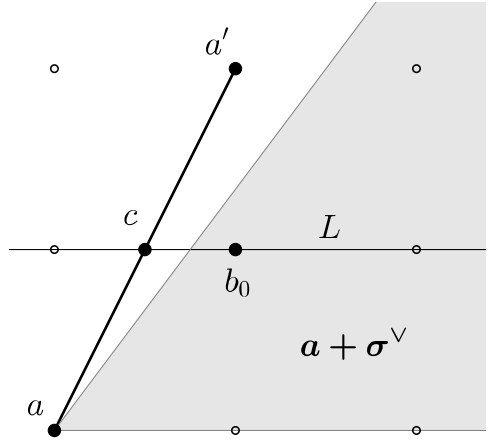}
				\caption{Neighborhood of $ b_0 $} \label{nbhofb0}
			\end{figure}
		\end{proof}		
	\end{Lem}
	
	Next, we will see a relation between $ \mathcal{P}_{n} $ and $ \mathcal{P}_{n+1} $.

	\begin{Lem}\label{Pn-lemma}
		Let $ n > 0 $ be an integer.
		\begin{enumerate}
			\item $ \sharp \mathcal{P}_n = n+3 $.
			\item For any distinct $ a,b \in \mathcal{P}_n $, we have $ b \notin a + \sigma_{\mathbb{Z}} $.
			\item If $ n $ is odd, then
			\begin{displaymath}
			p_{n-1} + (1,0) = p_{n} = p_{n+1}.
			\end{displaymath}
			If $ n $ is even, then
			\begin{displaymath}
			s_{n-1} + (3,4) = s_n = s_{n+1}.
			\end{displaymath}
			\item Consider $ \theta: \sigma^\vee \ni a \mapsto (1,1) + a \in \sigma^\vee $.	Then
			\begin{displaymath}
			\theta\left( q_n^i \right) = q_{n+1}^i,\ \theta\left( r_n^j\right) = r_{n+1}^j.
			\end{displaymath}
			Moreover if $ n $ is odd, then
			\begin{displaymath}
				\theta\left( s_n \right) = r_{n+1}^{\frac{n+1}{2}};
			\end{displaymath}
			if $ n $ is even, then 
			\begin{displaymath}
				\theta\left( p_n \right) = q_{n+1}^{\frac{n}{2}}.
			\end{displaymath}
			\item We have
			\begin{displaymath}
			\mathcal{P}_n \cap \mathcal{P}_{n+1} = \begin{cases}
			\left\lbrace \, p_n \, \right\rbrace & \text{if n is odd,}\\
			\left\lbrace \, s_n \, \right\rbrace & \text{if n is even}
			\end{cases}
			\end{displaymath}
			and 
			\begin{displaymath}
			\mathcal{P}_{n+1} = \theta\left( \mathcal{P}_n\setminus \mathcal{P}_{n+1}  \right) \sqcup \left\lbrace \, p_{n+1},s_{n+1} \, \right\rbrace.
			\end{displaymath}
			\item We have 
			\begin{displaymath}
			\mathcal{P}_{n} + \sigma_{\mathbb{Z}} = \left( \mathcal{P}_{n}\setminus\mathcal{P}_{n+1} \right) \sqcup \left(\mathcal{P}_{n+1} + \sigma_{\mathbb{Z}} \right).
			\end{displaymath}
		\end{enumerate}
	
		\begin{proof}
			(1) follows from direct calculations.
			
			(2) follows from Figure \ref{rn0+sigmaforodd} and Figure \ref{rn0+sigmaforeven}.
			
			(3) follows from direct calculations.
			
			(4) One can easily check that $ \theta\left( q_n^0 \right) = q_{n+1}^0 $ for any $ n $. Therefore we have
			\begin{displaymath}
			\theta\left( q_n^i \right) = \theta\left( q_n^0 - i\left( 1,2 \right) \right) = \theta\left( q_n^0 \right) - i\left( 1,2 \right) = q_{n+1}^0 - i\left( 1,2 \right) = q_{n+1}^i
			\end{displaymath}
			
			The other assertions follow from similar direct calculations with Lemma \ref{innerPn-lemma}.
			
			(5) By (3) and (4), we can describe the relation between $ \mathcal{P}_n $ and $ \mathcal{P}_{n+1} $ as in Figure \ref{Pn-relation-odd} and Figure \ref{Pn-relation-even}. In the figures, $ \theta $ shifts segments as follows; if $ n $ is odd, then
			\begin{displaymath}
			q_n^\frac{n-1}{2}q_n^0 \mapsto q_{n+1}^\frac{n-1}{2}q_{n+1}^0,\ r_n^0 s_n \mapsto r_{n+1}^0r_{n+1}^\frac{n}{2};
			\end{displaymath}
			if $ n $ is even, then
			\begin{displaymath}
			p_nq_n^0 \mapsto q_{n+1}^\frac{n}{2}q_{n+1}^0,\ r_n^0r_n^\frac{n}{2} \mapsto r_{n+1}^0r_{n+1}^\frac{n}{2}.
			\end{displaymath}
			Then one can easily check the assertions by the figures.
			
			\begin{figure}[htbp]
				\begin{tabular}{cc}
					\begin{minipage}{0.5\hsize}
						\centering
						\includegraphics[keepaspectratio, height=0.5\vsize]{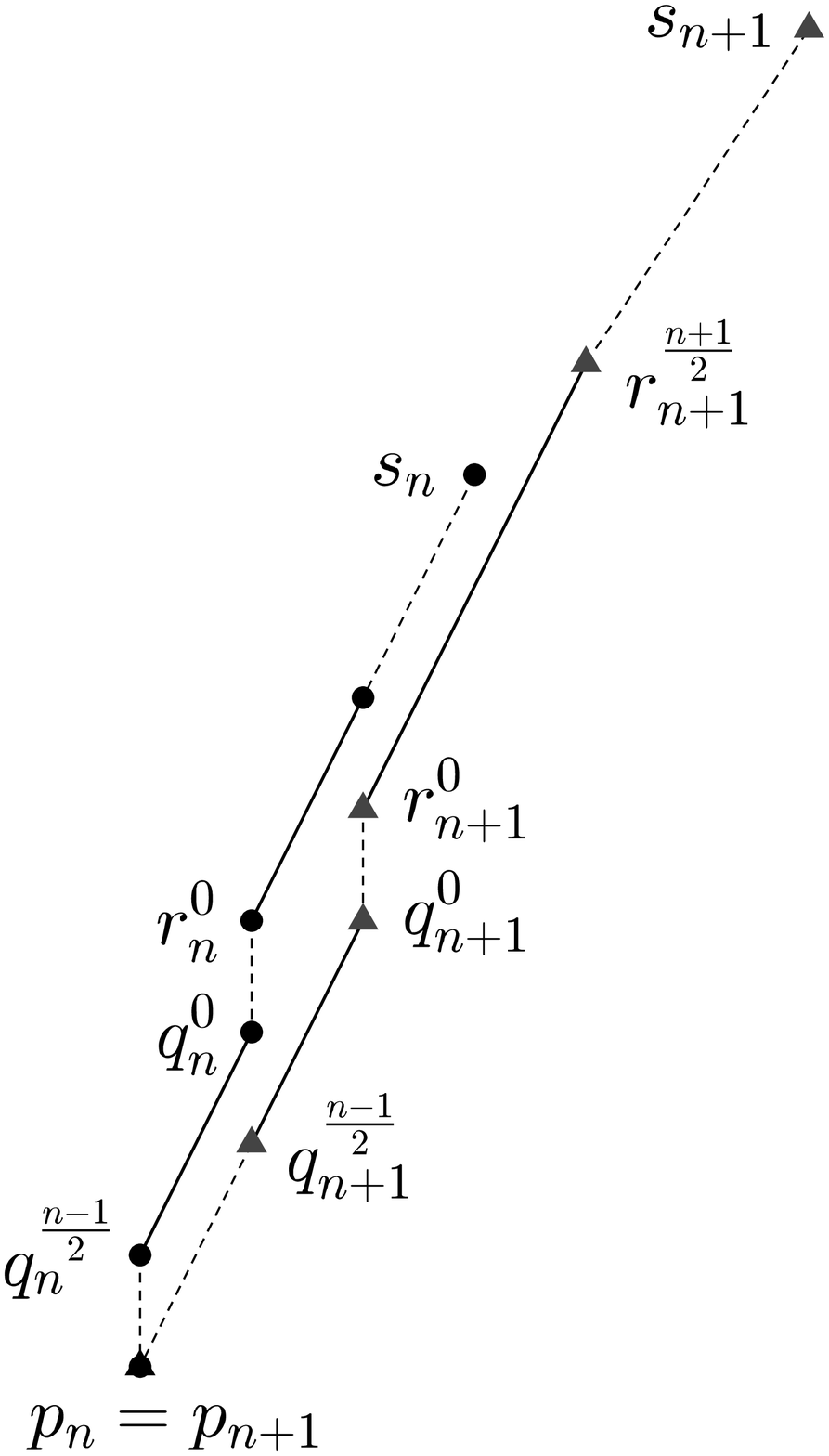}
						\caption{$ \mathcal{P}_n $ and $ \mathcal{P}_{n+1} $ for odd $ n $} \label{Pn-relation-odd}
					\end{minipage}
					\begin{minipage}{0.5\hsize}
						\centering
						\includegraphics[keepaspectratio, height=0.5\vsize]{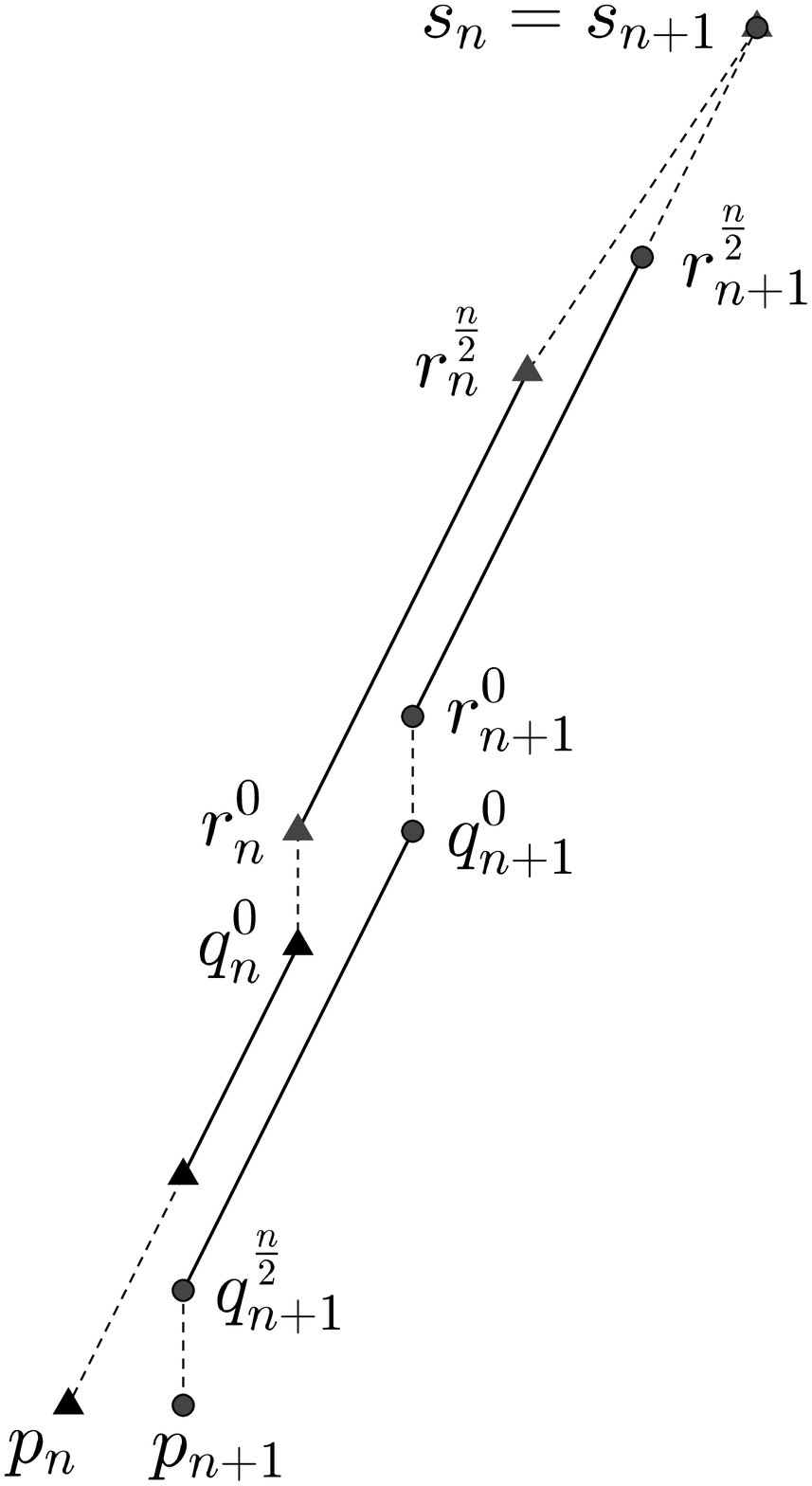}
						\caption{$ \mathcal{P}_n $ and $ \mathcal{P}_{n+1} $ for even $ n $} \label{Pn-relation-even}
					\end{minipage}
				\end{tabular}
			\end{figure}
			
			(6) By Figure \ref{Pn-relation-odd}, Figure \ref{Pn-relation-even} and Lemma \ref{Pn+sigma-lemma}, it is clear that $ \mathcal{P}_{n+1} \subset \mathcal{P}_{n} + \sigma_{\mathbb{Z}} $ and hence $ \mathcal{P}_{n+1} + \sigma_{\mathbb{Z}} \subset \mathcal{P}_{n} + \sigma_{\mathbb{Z}} $. Therefore it is enough to  see
			\begin{displaymath}
			\left(\mathcal{P}_{n} + \sigma_{\mathbb{Z}}\right) \setminus \left(\mathcal{P}_{n+1} + \sigma_{\mathbb{Z}}\right) = \mathcal{P}_n \setminus \mathcal{P}_{n+1}.
			\end{displaymath}
			
			Fix any $ a \in \mathcal{P}_n \setminus \mathcal{P}_{n+1} $. Then it is enough to check that
			\begin{displaymath}
			a \notin \mathcal{P}_{n+1} + \sigma_{\mathbb{Z}}\ \text{but}\ a + (1,0) , a + (3,4), a + (1,1) \in \mathcal{P}_{n+1} + \sigma_{\mathbb{Z}}
			\end{displaymath}
			because $ \sigma_{\mathbb{Z}} $ is generated by $ (1,0),(3,4),(1,1) $.
			
			$ a \notin \mathcal{P}_{n+1} + \sigma_{\mathbb{Z}} $ follows from Figure \ref{Pn-relation-odd} and Figure \ref{Pn-relation-even}.
			
			Let $ d $ be $ (1,0) $ or $ (3,4) $ or $ (1,1) $, and let us check $ a + d \in \mathcal{P}_{n+1} + \sigma_{\mathbb{Z}} $.
			
			The cases of $ d = (1,0),(1,1) $ are easily checked by Figure \ref{Pn-relation-odd} and Figure \ref{Pn-relation-even}. Let $ d = (3,4) $. Then
			\begin{displaymath}
			a + (3,4) = \theta(a) + (2,3)
			\end{displaymath}
			 and $ \theta(a) \in \mathcal{P}_{n+1} $. One can easily check that $ \alpha + (2,3) \in \mathcal{P}_{n+1} + \sigma_\mathbb{Z} $ for any $ \alpha \in \mathcal{P}_{n+1} $ by Figure \ref{Pn-odd} and Figure \ref{Pn-even} and Lemma \ref{Pn+sigma-lemma}. Thus $ a + d = \theta(a) + (2,3) \in \mathcal{P}_{n+1} + \sigma_\mathbb{Z} $.
		\end{proof}
	\end{Lem}

	\begin{Def}\label{Dn-definition}
		$ \mathcal{D}_n = \sigma_{\mathbb{Z}} \setminus \left( \mathcal{P}_n + \sigma_{\mathbb{Z}} \right) $. In $ S $, $ \mathcal{D}_n $ is the set of monomials not contained in the ideal $ \langle \, \mathcal{P}_n \, \rangle $. 
	\end{Def}

	\begin{Lem}\label{Dn-lemma} Let $ n > 0 $ be an integer.
		\begin{enumerate}
			\item $ \mathcal{D}_1 = \left\lbrace \, \left(0,0\right),\left(1,0\right),\left(1,1\right) \, \right\rbrace $ and $ \mathcal{D}_n = \mathcal{D}_{n-1} \sqcup \mathcal{P}_{n-1} \setminus \mathcal{P}_n $ for $ n \geq 2 $.
			\item $ \sharp \mathcal{D}_n = \frac{1}{2}(n+1)(n+2) $.
			\item $ (1,1) + \mathcal{D}_n \subset \mathcal{D}_{n+1} $.
			\item Let $ \Phi:\ \sigma_\mathbb{Z} \ni a \mapsto \left( -1,1 \right)\cdot a \in \mathbb{Z} $. If $ n $ is odd, then 
			\begin{displaymath}
			\Phi\left( \mathcal{D}_n \right) = \left\lbrace \, -\frac{n+1}{2}, -\left( \frac{n+1}{2}-1 \right), \ldots, -1, 0, 1, \ldots, \frac{n+1}{2}-1 \, \right\rbrace.
			\end{displaymath}
			If $ n $ is even, then
			\begin{displaymath}
			\Phi\left( \mathcal{D}_n \right) = \left\lbrace \, -\frac{n}{2}, -\left( \frac{n}{2}-1 \right), \ldots, -1, 0, 1, \ldots, \frac{n}{2}-1, \frac{n}{2} \, \right\rbrace
			\end{displaymath}
			and, in this case, $ s_{n-1} $ is the only member of $ \mathcal{D}_n $ whose image by $ \Phi $ is $ \frac{n}{2} $. 
			\item $ p_n $ is bigger than any monomial of $ \mathcal{D}_n $ with respect to $ \preceq $.
			\item Let $ n \geq 2 $. We define $ \Psi_n:\ \sigma_\mathbb{Z} \rightarrow \mathbb{R}_{\geq 0} $ by $ a \mapsto l_n \cdot a $ for
			$$ l_n :=
			\begin{cases}
			(2n-2,-n+2) & \text{if}\ n\ \text{is odd},\\
			(2n,-n+1) & \text{if}\ n\ \text{is even}.\\
			\end{cases}
			$$ Then if $ n $ is odd, we have
			\begin{displaymath}
			\max \Psi_n \left( \mathcal{D}_n \right) = \Psi_n\left(r_{n-1}^{\frac{n-1}{2}}\right)
			\end{displaymath}
			and this is also equal to
			\begin{displaymath}
			\min \Psi_n \left( \mathcal{P}_n \right) = \Psi_n\left(q_n^{\frac{n-1}{2}}\right).
			\end{displaymath}
			If $ n $ is even, we have
			\begin{displaymath}
			\max \Psi_n \left( \mathcal{D}_n \right) = \Psi_n\left(s_{n-1}\right)
			\end{displaymath}
			and this is also equal to
			\begin{displaymath}
			\min \Psi_n \left( \mathcal{P}_n \right) = \Psi_n\left(p_n\right).
			\end{displaymath}
		\end{enumerate}
	
		\begin{proof} (1) The assertion for $ n=1 $ follows from Figure \ref{P1}.
			
			\begin{figure}[htbp]
				\centering
				\includegraphics[keepaspectratio, height=0.7\hsize]{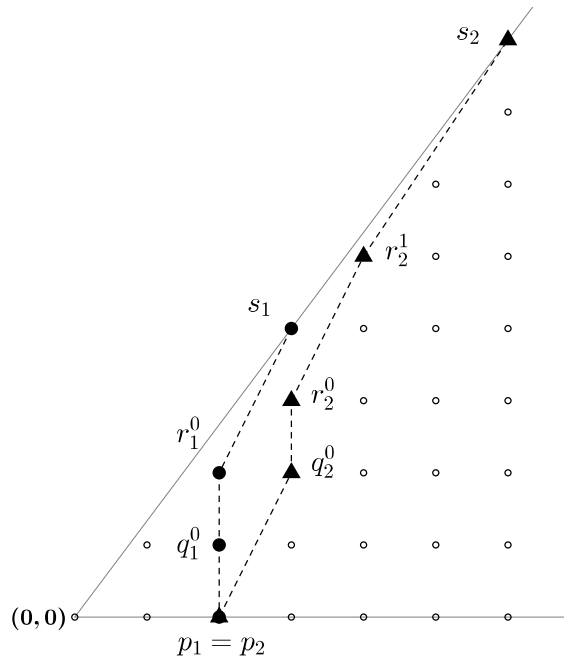}
				\caption{$ \mathcal{P}_1 $ and $ \mathcal{P}_2 $}\label{P1}
			\end{figure}
			
			Let $ n \geq 2 $. Then Lemma \ref{Pn-lemma} (6) implies the following equations;
			\begin{align*}
			\mathcal{D}_n & = \sigma_{\mathbb{Z}} \setminus \left( \mathcal{P}_n + \sigma_{\mathbb{Z}} \right)\\
			& = \left(\sigma_{\mathbb{Z}} \setminus \left( \mathcal{P}_{n-1} + \sigma_{\mathbb{Z}} \right)\right) \sqcup \left(\mathcal{P}_{n-1} + \sigma_{\mathbb{Z}}\right) \setminus \left(\mathcal{P}_{n} + \sigma_{\mathbb{Z}} \right)\\
			& = \mathcal{D}_{n-1} \sqcup \mathcal{P}_{n-1} \setminus \mathcal{P}_n.
			\end{align*}
			
		Hence (1) holds.
		
		The proofs of (2)-(4) are by induction on $ n $. (2) and (3) for $ n=1 $ follow from Figure \ref{P1}, and (4) for $ n=1 $ follows from direct calculations.
			
		(2) Assume $ n \geq 2 $. Now $ \sharp \mathcal{P}_{n-1} \setminus \mathcal{P}_n = n+1 $ by Lemma \ref{Pn-lemma} (1) and Lemma \ref{Pn-lemma} (5). Hence (1) implies
		\begin{displaymath}
			\sharp \mathcal{D}_n = \sharp \mathcal{D}_{n-1} + \sharp\mathcal{P}_{n-1} \setminus \mathcal{P}_n = \frac{1}{2}n(n+1) + \left(n+1\right) = \frac{1}{2}(n+1)(n+2).
		\end{displaymath}
			
		(3) Assume $ n \geq 2 $. By (1) we have
		\begin{displaymath}
		(1,1) + \mathcal{D}_n = \left((1,1) + \mathcal{D}_{n-1}\right) \cup \left((1,1) + \left(\mathcal{P}_{n-1} \setminus \mathcal{P}_n\right)\right)
		\end{displaymath}
		Now (1) also shows $ \mathcal{D}_n \subset \mathcal{D}_{n+1} $. Thus, by the induction hypothesis, we have
		\begin{displaymath}
			(1,1) + \mathcal{D}_{n-1} \subset \mathcal{D}_n \subset \mathcal{D}_{n+1}.
		\end{displaymath}
		 Moreover Lemma \ref{Pn-lemma} (5) shows that
		\begin{displaymath}
			(1,1) + \left(\mathcal{P}_{n-1} \setminus \mathcal{P}_n\right) = \theta\left( \mathcal{P}_{n-1} \setminus \mathcal{P}_n \right) \subset \mathcal{P}_n \setminus \mathcal{P}_{n+1} \subset \mathcal{D}_{n+1}.
		\end{displaymath}
			Therefore $ (1,1) + \mathcal{D}_n \subset \mathcal{D}_{n+1} $.
			
		(4) Assume $ n \geq 2 $. By (1) we have $ \Phi\left( \mathcal{D}_n \right) = \Phi\left( \mathcal{D}_{n-1} \right) \cup \Phi \left( \mathcal{P}_{n-1} \setminus \mathcal{P}_n \right) $.
			
		Let $ n $ be odd. By the induction hypothesis, $ \Phi\left( \mathcal{D}_{n-1} \right) $ consists of
		\begin{displaymath}
		-\left(\frac{n+1}{2}-1\right), -\left( \frac{n+1}{2}-2 \right), \ldots, -1, 0, 1, \ldots, \frac{n+1}{2}-2, \frac{n+1}{2}-1.
		\end{displaymath}
		Moreover $ \Phi\left( \mathcal{P}_{n-1} \setminus \mathcal{P}_n \right) $ consists of
		\begin{displaymath}
		\Phi\left( p_{n-1} \right) = -\frac{n+1}{2},\ \Phi\left( q_{n-1}^i \right) = -1 - i,\ \Phi\left( r_{n-1}^j \right) = j
		\end{displaymath}
		where $ 0 \leq i \leq \frac{n-3}{2},\ 0 \leq j \leq \frac{n-1}{2} $. Hence the assertion holds for odd $ n $.
		
		Let $ n $ be even. By the induction hypothesis, $ \Phi\left( \mathcal{D}_{n-1} \right) $ consists of
		\begin{displaymath}
		-\frac{n}{2}, -\left( \frac{n}{2}-1 \right), \ldots, -1, 0, 1, \ldots, \frac{n}{2}-1.
		\end{displaymath}
		Moreover $ \Phi\left( \mathcal{P}_{n-1} \setminus \mathcal{P}_n \right) $ consists of 
		\begin{displaymath}
		\Phi\left( q_{n-1}^i \right) = -1 - i,\ \Phi\left( r_{n-1}^j \right) = j,\ \Phi\left(s_{n-1}\right) = \frac{n}{2}
		\end{displaymath}
		where $ 0 \leq i,j \leq \frac{n-2}{2} $. Hence the assertion holds.
		
		(5) By (1) we can describe $ \mathcal{D}_n $ as in Figure \ref{Dn-odd} and Figure \ref{Dn-even}; $ \mathcal{D}_n $ consists of lattice points of the shadow area with border.
		
		\begin{figure}[htbp]
			\begin{tabular}{cc}
				\begin{minipage}{0.5\hsize}
					\centering
					\includegraphics[keepaspectratio, width=\hsize]{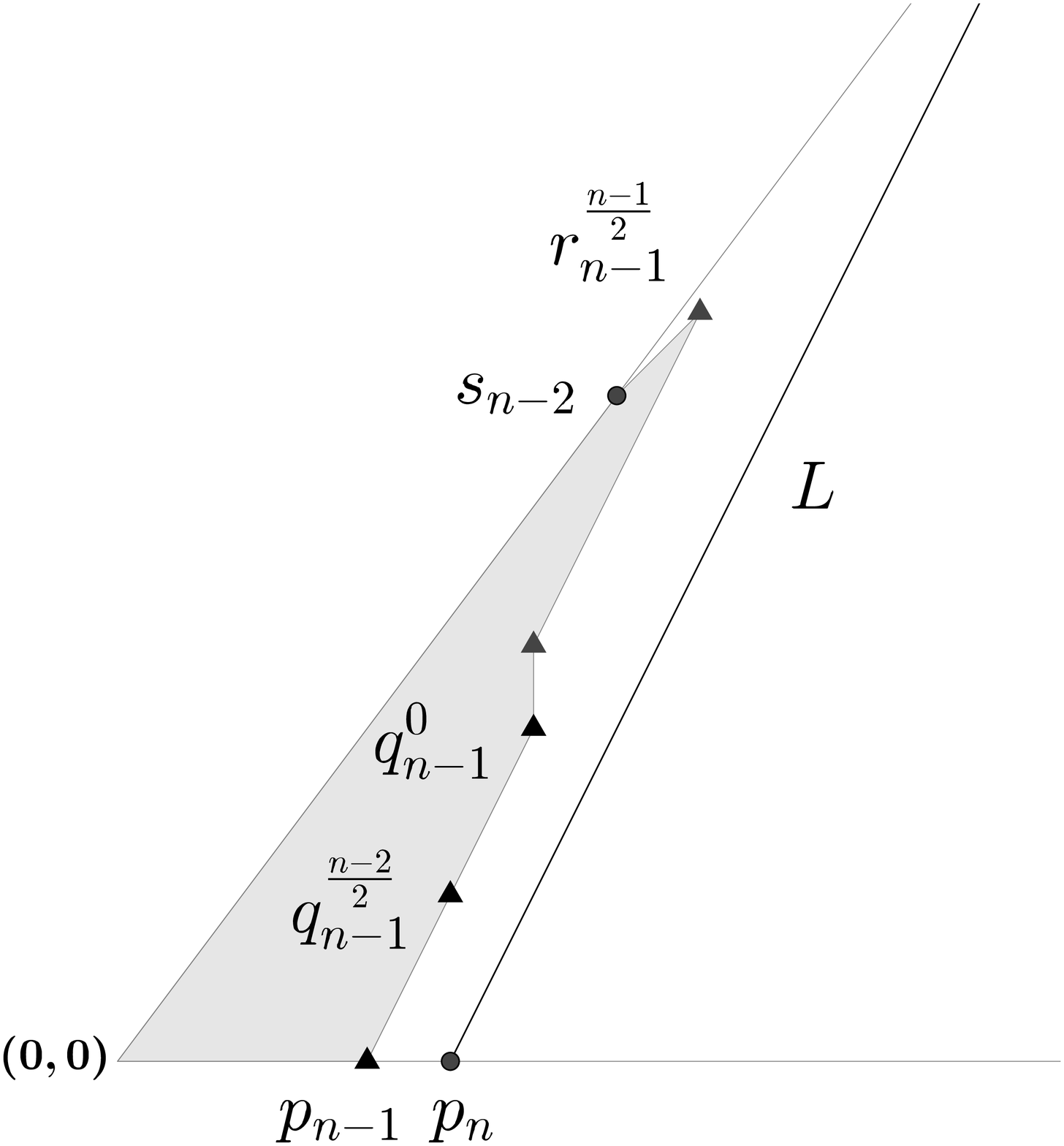}
					\caption{$ \mathcal{D}_n $ for odd $ n $} \label{Dn-odd}
				\end{minipage}
				\begin{minipage}{0.5\hsize}
					\centering
					\includegraphics[keepaspectratio, width=\hsize]{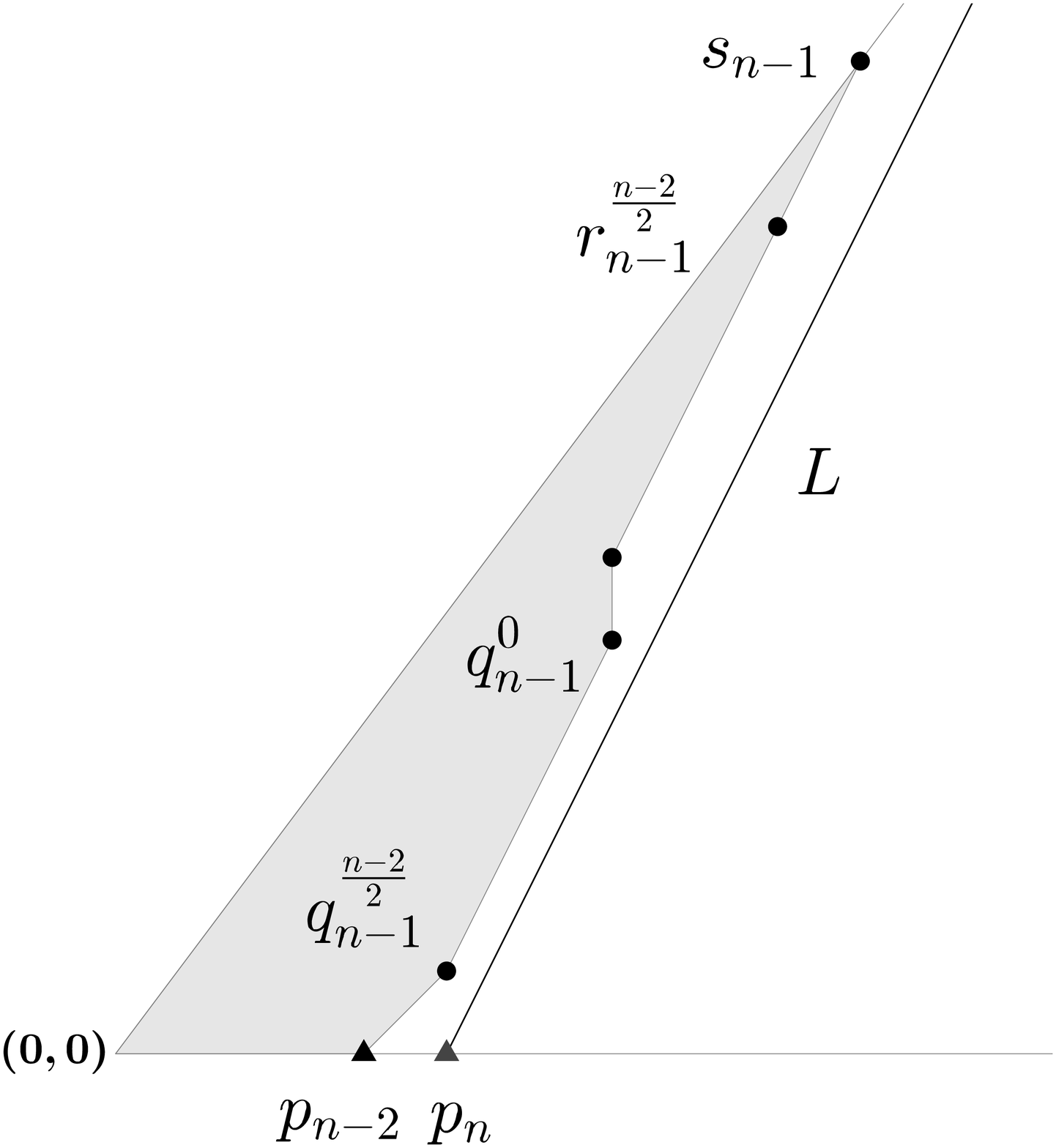}
					\caption{$ \mathcal{D}_n $ for even $ n $} \label{Dn-even}
				\end{minipage}
			\end{tabular}
		\end{figure}
		
		Recall the definition of $ \preceq $ (Definition \ref{<-definition}). Then we only have to see 
		\begin{displaymath}
			(2,-1)\cdot p_n > (2,-1)\cdot a\ \text{for all}\ a \in \mathcal{D}_n.
		\end{displaymath}
		Let $ L $ be the line passing through $ p_n $ whose normal vector $ (2,-1) $ (as in Figure \ref{Dn-odd} and Figure \ref{Dn-even}). Then $ \mathcal{D}_n $ is contained in the half plane whose border $ L $. This implies the assertion. 
		
		(6) Let $ a \in \mathcal{D}_n $ (resp.\ $ \mathcal{P}_n $). Then $ \Psi_n $ on $ \mathcal{D}_n $ (resp.\ $ \mathcal{P}_n $) attains the maximum (resp.\ minimum) value at $ a $ if and only if the folloing holds; let $ L' $ be the line passing through $ a $ whose normal vector $ l_n $, then $ \mathcal{D}_n $ (resp.\ $ \mathcal{P}_n $) is contained in the closed half plane $ H $ whose border $ L' $ with $ a - l_n \in H $ (resp.\ $ a + l_n \in H $).
		
		Let $ n $ be odd. Now the slope of $ L' $ is $ \frac{2n-2}{n-2} > 2 $. Thus, by Figure \ref{Dn-odd}, one can easily see
		\begin{displaymath}
			\max \Psi_n\left(\mathcal{D}_n\right) = \max \left\lbrace \, \Psi_n\left(r_{n-1}^{\frac{n-1}{2}}\right),\Psi_n\left(q_{n-1}^0\right) \, \right\rbrace.
		\end{displaymath}
		By direct calculations, we have
		\begin{displaymath}
			\max \Psi_n\left(\mathcal{D}_n\right) = \Psi_n\left(r_{n-1}^{\frac{n-1}{2}}\right) = (n-1)(n+2)+1.
		\end{displaymath}
		
		On the other hand, by Figure \ref{Pn-odd}, one can easily see
		\begin{displaymath}
			\min \Psi_n\left(\mathcal{P}_n\right) = \min \left\lbrace \, \Psi_n\left(q_n^\frac{n-1}{2}\right),\Psi_n\left(r_n^0\right) \, \right\rbrace.
		\end{displaymath}		
		By direct calculations, we have
		\begin{displaymath}
		\min \Psi_n\left(\mathcal{P}_n\right) = \Psi_n\left(q_n^\frac{n-1}{2}\right) = (n-1)(n+2) + 1
		\end{displaymath}
		and the assertion holds for odd $ n $.
		
		Let $ n $ be even. Now the slope of $ L' $ is $ \frac{2n}{n-1} >2 $. Thus, by Figure \ref{Dn-even}, one can easily see
		\begin{displaymath}
		\max \Psi_n\left(\mathcal{D}_n\right) = \max \left\lbrace \, \Psi_n\left(s_{n-1}\right),\Psi_n\left(q_{n-1}^0\right) \, \right\rbrace.
		\end{displaymath}
		By direct calculations, we have
		\begin{displaymath}
		\max \Psi_n\left(\mathcal{D}_n\right) = \Psi_n\left(r_{n-1}^{\frac{n-1}{2}}\right) = n(n+2).
		\end{displaymath}
		
		On the other hand, by Figure \ref{Pn-even}, one can easily see
		\begin{displaymath}
		\min \Psi_n\left(\mathcal{P}_n\right) = \min \left\lbrace \, \Psi_n\left(p_n\right),\Psi_n\left(r_n^0\right) \, \right\rbrace.
		\end{displaymath}		
		By direct calculations, we have
		\begin{displaymath}
		\min \Psi_n\left(\mathcal{P}_n\right) = \Psi_n\left(p_n\right) = n(n+2)
		\end{displaymath}
		and the assertion holds for even $ n $ also.
		\end{proof}
	\end{Lem}

	\subsection{Proof of $ \mathbb{M}_n = \mathcal{P}_n $}
	
	It is clear that $ \mathbb{M}_n $ generates $ \initial_{\preceq}\left(J_n\right) $. The key of the proof for $ \mathbb{M}_n = \mathcal{P}_n $ is to see that $ \initial_{\preceq}\left( J_n \right) $ is also generated by $ \mathcal{P}_n $. We prepare some lemmas for the proof.

	\begin{Lem}[\cite{Duarte2}, Appendix A, Proposition A.2.1]\label{Macaulay'sTh}
		For any ideal $ I $ of $ S $, the monomials of $ S $ not contained in $ \initial_{\preceq}(I) $ form a $ \mathbb{C} $-basis of $ S/I $. Therefore we have
		\begin{enumerate}
			\item $ \dim_\mathbb{C} S/J_n = \dim_\mathbb{C} S/\initial_{\preceq}\left(J_n\right) $,
			\item $ \dim_\mathbb{C} S/\langle \, \mathcal{P}_n \, \rangle = \sharp \mathcal{D}_n $.
		\end{enumerate}
	\end{Lem}

	\begin{Lem}\label{QuotDim-lemma}
		
		\begin{enumerate}
			\item $ \dim_\mathbb{C} S / \initial_{\preceq}\left(J_n\right) = \frac{1}{2}\left( n+1 \right)\left( n+2 \right) = \dim_\mathbb{C} S / \langle \, \mathcal{P}_n \, \rangle $.
			\item $ \left( J_n : uv-1 \right)_S = J_{n-1} $.
			\item $ \dim_\mathbb{C} \initial_{\preceq}\left(J_{n-1}\right) / \initial_{\preceq}\left(J_{n}\right) = n+1 $. Furthermore if a set of monomials generates $ \initial_{\preceq}\left(J_{n-1}\right) $ as an ideal, then the set generates $ \initial_{\preceq}\left(J_{n-1}\right) / \initial_{\preceq}(J_{n}) $ as a vector space over $ \mathbb{C} $.
		\end{enumerate}	
		
		\begin{proof}
			(1) By Lemma \ref{Macaulay'sTh} (1), we can consider $ \dim_\mathbb{C} S / J_n $ instead of $ \dim_\mathbb{C} S / \initial_{\preceq}\left(J_n\right) $.
			
			Let $ J_0 := \langle \, u-1,u^3v^4-1,uv-1 \, \rangle $. Then $ S_{J_0} $ is a regular local ring because $ J_0 $ is the maximal ideal corresponding to the regular point $ \left(1,1,1\right) $ of $ X = (z^4-xy = 0) $. Moreover $ J_0S_{J_0} = \langle \, u-1, uv-1 \, \rangle $ since the following equation holds;
			\begin{equation*}
			u^3v^4 - 1 = \left(u^3v^3 + u^2v^2 + uv + 1\right)\left(uv - 1\right) - u^3v^4\left(u-1\right).
			\end{equation*}
			
			Now consider $ \text{gr}_{J_0}(S) = \bigoplus_{\nu=0}^\infty J_0^\nu / J_0^{\nu+1} $. Then we obtain an isomorphism of graded rings
			\begin{align*}
			&\mathbb{C}[x_1,x_2] \overset{\cong}{\longrightarrow} \text{gr}_{J_0S_{J_0}}(S_{J_0}) \cong \text{gr}_{J_0}(S);\\
			 & x_1 \mapsto \left[ u-1 \ \text{mod}\, J_0^2 \right],\ x_2 \mapsto \left[uv-1\ \text{mod}\, J_0^2\right].
			\end{align*}
			
			Hence 
			\begin{displaymath}
			\dim_\mathbb{C} S / J_n = \dim_\mathbb{C} \mathbb{C}[x_1,x_2] / \langle \, x_1,x_2 \, \rangle^{n+1} = \frac{1}{2}\left( n+1 \right)\left( n+2 \right).
			\end{displaymath}
			The last equation follows from Lemma \ref{Dn-lemma} (2) and Lemma \ref{Macaulay'sTh} (2).
			
			(2) It is clear that $ \left( J_n : uv-1 \right)_S \supset J_{n-1} $. Fix any $ f \in \left( J_n : uv-1 \right)_S $ and let us show $ f \in J_{n-1} $. 
			
			Consider $ \text{gr}_{J_0}(S) $ again. Suppose that $ f \in J_0^i $ for $ 0 \leq i \leq n-1 $. Then we have 
			\begin{displaymath}
			\left[uv-1\ \text{mod}\, J_0^2 \right] \cdot \left[ f\ \text{mod}\, J_0^{i+1} \right] = \left[ (uv-1)f\ \text{mod}\, J_0^{i+2} \right] = 0
			\end{displaymath}
			because $  (uv-1)f \in  J_n \subset J_{i+1} = J_0^{i+2} $. However $ \text{gr}_{J_0}(S) $ is an integral domain as above, hence $ f = 0\ \text{mod}\, J_0^{i+1} $. Therefore $ f \in J_0^{i+1} $, so $ f \in J_0^{n} = J_{n-1} $.
			
			(3) By (1), we have
			\begin{align*}
			\dim_\mathbb{C} \initial_{\preceq}\left(J_{n-1}\right) / \initial_{\preceq}(J_{n}) & = \dim_\mathbb{C} S /  \initial_{\preceq}(J_{n})  - \dim_\mathbb{C} S /  \initial_{\preceq}(J_{n-1})\\
			& = \frac{1}{2}\left( n+1 \right)\left( n+2 \right) - \frac{1}{2}n\left( n+1 \right)\\
			& = n+1.
			\end{align*}
			
			For the last assertion, let $ \left\lbrace \, m_1,\ldots,m_r \, \right\rbrace $ be any set of monomials generating $ \initial_{\preceq}\left(J_{n-1}\right) $ as an ideal. The vector space $ \initial_{\preceq}\left(J_{n-1}\right) / \initial_{\preceq}(J_{n}) $ is generated by monomials of $ \initial_{\preceq}\left(J_{n-1}\right) $. Let $ m $ be any monomial of $ \initial_{\preceq}\left(J_{n-1}\right) $. Then $ m $ is divisible by some $ m_i $. If $ m\neq m_i $, then there exists $ u^av^b \in \left\lbrace \, u,u^3v^4,uv \, \right\rbrace $ such that $ m $ is divisible by $ m_i(u^av^b) $. However one can fined $ f \in J_{n-1} $ with $ \lm{f} = m_i $ and obtain $ g:=(u^av^b-1)f \in J_{n} $. Then $ m_i(u^av^b) = \lm{g} \in \initial_{\preceq}(J_{n}) $ and hence $ m = 0 $ in $ \initial_{\preceq}\left(J_{n-1}\right) / \initial_{\preceq}(J_{n}) $. Therefore $ \initial_{\preceq}\left(J_{n-1}\right) / \initial_{\preceq}(J_{n}) $ is generated by $ \left\lbrace \, m_1,\ldots,m_r \, \right\rbrace $ as a vector space.
		\end{proof}
	\end{Lem}

	The following proposition determines $ \mathbb{G}_1 $.
	
	\begin{Prop}\label{n=1}
		The reduced Gr\"{o}bner basis of $ J_1 $ with respect to $ \preceq $ consists of the following polynomials;
		\begin{align*}
		& \underline{u^3v^4}+u-4uv+2,\\
		& \underline{u^2v^2}-2uv+1,\\
		& \underline{u^2v} - u -uv +1,\\
		& \underline{u^2} - 2u +1,
		\end{align*}
		where the underlined monomials are the leading terms with respect to $ \preceq $. Therefore $ \mathbb{M}_1 $ coincides with $ \mathcal{P}_1 = \left\lbrace \, \left(3,4\right), \left(2,2\right), \left(2,1\right), \left(2,0\right) \, \right\rbrace $.
		
		\begin{proof}
			First we will see that the polynomials are contained in $ J_1 $ and their leading terms are as asserted.
			
			Let $ g_1 := u^3v^4+u-4uv+2 $. Then
			\begin{displaymath}
			g_1 = \left(\left(uv\right)^2 + 2uv +3\right)\left(uv-1\right)^2 - \left(u-1\right)\left(u^3v^4-1\right) \in J_1.
			\end{displaymath}
			Moreover
			\begin{align*}
			&\left(2,-1\right)\cdot\left(3,4\right) = 2 \geq \left(2,-1\right)\cdot\left(1,0\right) = 2 > \left(2,-1\right)\cdot\left(1,1\right) = 1,\\
			&\left(1,1\right)\cdot\left(3,4\right) = 7 > \left(1,1\right)\cdot\left(1,0\right) = 1,
			\end{align*}
			and hence $ \lm{g_1} = u^3v^4 $. It is easy to see for the other polynomials;
			\begin{align*}
			& g_2 := u^2v^2-2uv+1 = \left(uv-1\right)^2 \in J_1,\\
			& g_3 := u^2v - u -uv +1 = \left(u-1\right)\left(uv-1\right) \in J_1,\\
			& g_4 := u^2 - 2u +1 = \left(u-1\right)^2 \in J_1.
			\end{align*}
			
			One can easily check that there exists no monomial of $ \supp{g_i} $ which is divisible by $ \lm{g_j} $ for $ j \neq i $. Therefore we only have to show
			\begin{displaymath}
			\initial_{\preceq}\left(J_1\right) = \langle \, u^3v^4, u^2v^2, u^2v, u^2 \, \rangle.
			\end{displaymath}
			We have $ \initial_{\preceq}\left(J_1\right) \supset \mathcal{P}_1 $ and hence it is enough to see $ \dim_\mathbb{C} S / \initial_{\preceq}\left(J_1\right) = \dim_\mathbb{C} S / \langle \, \mathcal{P}_1 \, \rangle $.	Thus Lemma \ref{QuotDim-lemma} (1) completes the proof.
		\end{proof}
	\end{Prop}

	The cases of higher $ n > 0 $ need the next lemma.
	
	\begin{Lem}\label{phi-lemma}
		Consider the morphism of $ \mathbb{C} $-algebras $ \mathbb{C}[u,v] \rightarrow \mathbb{C}[\lambda,\lambda^{-1}] $ given by $ u \mapsto \lambda^{-1}, v \mapsto \lambda $. By restriction to $ S \subset \mathbb{C}[u,v] $, we obtain
		\begin{displaymath}
			\phi: S \rightarrow \mathbb{C}[\lambda^\pm];\ u \mapsto \lambda^{-1},\ u^3v^4 \mapsto \lambda,\ uv \mapsto 1.
		\end{displaymath}
		Then
		\begin{enumerate}
			\item  $ \phi $ is a surjection and $ \ker(\phi) = \langle \, uv-1 \, \rangle $.
			
			\item Let $ n>0 $ be even. Then
			\begin{displaymath}
			B := \left\lbrace \, \lambda^{-\frac{n}{2}}, \lambda^{-\left( \frac{n}{2}-1 \right)}, \ldots, \lambda^{-1}, 1, \lambda, \ldots, \lambda^{\frac{n}{2}-1}, \lambda^{\frac{n}{2}} \, \right\rbrace \subset \mathbb{C}[\lambda^{\pm}]
			\end{displaymath}
			forms a $ \mathbb{C} $-basis of $ \mathbb{C}[\lambda^{\pm}] / \phi(J_n) $.
			
			\item Let $ n>0 $ be even. Let $ f \in S $ satisfy $ \lm{f} \in \mathcal{P}_{n-1} \setminus \mathcal{P}_n $. Then
			\begin{displaymath}
			\supp{\phi(f)} \subset B.
			\end{displaymath}
			Therefore, if such $ f $ is in $ J_n $, then $ \phi\left(f\right) = 0 $ by (2).
			
		\end{enumerate}
		
		\begin{proof}
			(1) $ \phi $ is obviously a surjection. Recall $ F: \mathbb{C}[x,y,z] \twoheadrightarrow S $ in Notation \ref{Settings} (3). Then one can easily see $ \ker \left( F\circ\phi \right) = \langle \, xy-1,z-1 \, \rangle $. Therefore
			\begin{displaymath}
			\ker \phi = F\left( \langle \, xy-1,z-1 \, \rangle \right) = \langle \, u^4v^4 -1, uv-1 \, \rangle = \langle \, uv-1 \, \rangle
			\end{displaymath}
			because  $ u^4v^4 -1 = \left(u^3v^3 + u^2v^2 + uv + 1\right)\left(uv - 1\right) $.
			
			(2) One can easily check
			\begin{displaymath}
			\phi\left( u - 1 \right) = -\lambda^{-1}\phi\left( u^3v^4-1 \right),\ \phi\left( uv-1 \right) = 0.
			\end{displaymath}
			Thus we have $ \phi(J_n) = \langle \, \phi\left(u^3v^4-1\right) \, \rangle^{n+1} = \langle\,  \lambda-1 \, \rangle^{n+1} $.
			
			It is clear that $ 1,\lambda,\ldots,\lambda^n $ form a $ \mathbb{C} $-basis of $ \mathbb{C}[\lambda] / \langle \lambda-1 \rangle^{n+1} $ and
			\begin{displaymath}
			\mathbb{C}[\lambda] / \langle \, \lambda-1 \, \rangle^{n+1} = \mathbb{C}[\lambda^\pm] / \langle \, \lambda-1 \, \rangle^{n+1} = \mathbb{C}[\lambda^{\pm}] / \phi(J_n).
			\end{displaymath}
			Now $ \lambda $ is a unit element of the quotient rings. Hence, by multiplying $ \lambda^{-\frac{n}{2}} $, we obtain $ B = \lambda^{-\frac{n}{2}} \cdot \left\lbrace \, 1,\lambda,\ldots,\lambda^n \, \right\rbrace  $ as a $ \mathbb{C} $ basis of $ \mathbb{C}[\lambda^{\pm}] / \phi(J_n) $.
			
			(3) Fix any $ m \in \supp{f} $ and let us show $ \phi\left(m\right) \in B $.
			
			Suppose $ \phi(m) \notin B $. If $ \phi(m)=\lambda^d $ for $ d \geq \frac{n}{2}+1 $, then $ m $ is divisible by $ \left( u^3v^4 \right)^d $. Indeed $ m $ can be written as $ m = u^a\left(u^3v^4\right)^b\left(uv\right)^c $ for some $ a,b,c \geq 0 $, and we have $ \phi\left( u \right) = \lambda^{-1} $, $ \phi\left( u^3v^4 \right) = \lambda $ and $ \phi\left( uv \right) = 1 $. Therefore $ b \geq d $.  
			
			However this leads to the contradiction $ \lm{f} \prec m $. Indeed $ \lm{f} \in \mathcal{P}_{n-1}\setminus\mathcal{P}_n = \mathcal{P}_{n-1}\setminus\left\lbrace \, p_{n-1} \, \right\rbrace $, and hence 
			$$
			\lm{f} \cdot (2,-1) =
			\begin{cases}
			n+1 & \text{if }\lm{f} = q_{n-1}^i,\\
			n & \text{if }\lm{f} = r_{n-1}^i\ \text{or } s_{n-1}.
			\end{cases}
			$$
			On the other hand, $ m \cdot (2,-1) \geq (3d,4d)\cdot(2,-1) = 2d \geq n+2 $. Therefore $ \lm{f} \prec m $ by the definition of $ \preceq $. This contradicts to $ m \in \supp{f} $.
			
			If $ \phi(m)=\lambda^{-d} $ for $ d \geq \frac{n}{2}+1 $, then $ m $ is divisible by $ u^d $. Now $ (d,0)\cdot(2,-1)=2d \geq n+2 $ and the contradiction $ \lm{f} \prec m $ is also induced.
			
			Therefore $ \phi\left(m\right) \in B $.
		\end{proof}
	\end{Lem}

	The following proposition is the first consequence of above arguments.

	\begin{Prop}\label{lm-proposition}
		Let $ n >0 $ be an integer. Then $ \mathbb{M}_n $ coincides with $ \mathcal{P}_n $.
		
		\begin{proof}
			By induction on $ n $. The case of $ n=1 $ has already been done.
			
			Let $ n \geq 2 $. The arguments will go as follows. We will show $ \mathcal{P}_n \subset \initial_{\preceq}\left(J_n\right) $. Then we can conclude that $ \initial_{\preceq}\left(J_n\right) $ is generated by the monomials of  $ \mathcal{P}_n $by Lemma \ref{QuotDim-lemma} (1). Therefore, we can conclude $ \mathbb{M}_n = \mathcal{P}_n $ by Lemma \ref{Pn-lemma} (2).
			
			Let us show $ \mathcal{P}_n \subset \initial_{\preceq}\left(J_n\right) $. Lemma \ref{Pn-lemma} (5) shows that
			\begin{displaymath}
			\mathcal{P}_{n} = \begin{cases}
			\theta\left( \mathcal{P}_{n-1}\setminus \left\lbrace \, s_{n-1} \, \right\rbrace  \right) \sqcup \left\lbrace \, p_{n},s_{n} \, \right\rbrace & \text{if n is odd,}\\
			\theta\left( \mathcal{P}_{n-1}\setminus \left\lbrace \, p_{n-1} \, \right\rbrace  \right) \sqcup \left\lbrace \, p_{n},s_{n} \, \right\rbrace & \text{if n is even}.\\
			\end{cases}
			\end{displaymath}
			
			To begin with, let $ n $ be odd.
			
			To see $ \theta\left( \mathcal{P}_{n-1}\setminus \left\lbrace \, s_{n-1} \, \right\rbrace  \right) \subset \initial_{\preceq}\left(J_n\right) $, fix any $ \alpha \in \mathcal{P}_{n-1} $ and let us show $ \theta\left( \alpha \right) \in \initial_{\preceq}\left( J_n \right) $. By the induction hypothesis, we have $ \alpha \in \mathbb{M}_{n-1} $. Thus there exists $ \left(f ,\alpha\right) \in \mathbb{G}_{n-1} $. Now $ f \in J_{n-1} $ and hence $ (uv-1)f \in J_n $. Therefore $ \theta\left( \alpha \right) = \lm{(uv-1)f} \in \initial_{\preceq}\left(J_n\right) $.
			
			To see $ s_n \in \initial_{\preceq}\left(J_n\right) $, let us remark that $ J_1 $ has $ g := u^3v^4+u-4uv+2 $ with $ \lm{g} = (3,4) $ as in Proposition \ref{n=1}. Therefore $ g^{\frac{n+1}{2}} \in \left(J_1\right)^{\frac{n+1}{2}} = \left(J_0\right)^{n+1} = J_{n} $, and hence $ s_n = \frac{n+1}{2}\left( 3,4 \right) = \lm{g^{\frac{n+1}{2}}} \in \initial_{\preceq}\left(J_n\right) $.
			
			Furthermore $ p_n \in \initial_{\preceq}\left(J_n\right) $. Indeed, by the induction hypothesis, $ p_{n-1} \in \mathbb{M}_{n-1} $. Hence there exists $ \left( h, p_{n-1} \right) \in \mathbb{G}_{n-1} $. Then $ (u-1)h \in J_n $ and hence $ p_n = p_{n-1} + (1,0) = \lm{(u-1)h} \in \initial_{\preceq}\left(J_n\right) $.
			
			Therefore $ \mathcal{P}_n \subset \initial_{\preceq}\left(J_n\right) $. 
			
			Next, let $ n $ be even.
			
			One can see $ \theta\left( \mathcal{P}_{n-1}\setminus \left\lbrace \, p_{n-1} \, \right\rbrace  \right) \subset \initial_{\preceq}\left(J_n\right) $ by arguments similar to the above. Moreover $ s_n \in \initial_{\preceq}\left(J_{n}\right) $. Indeed, by the induction hypothesis, $ s_{n-1} \in \mathbb{M}_{n-1} $. Hence one can find $ \left(g, s_{n-1} \right) \in \mathbb{G}_{n-1} $. Then $ (u^3v^4-1)g \in J_n $ and $ s_n = s_{n-1} + (3,4) = \lm{(u^3v^4-1)g} \in \initial_{\preceq}\left(J_n\right) $.
			
			It is slightly hard to see $ p_n \in \initial_{\preceq}\left(J_n\right) $ as follows.
			
			Lemma \ref{QuotDim-lemma} (3) shows that $ \initial_{\preceq}\left(J_{n-1}\right) / \initial_{\preceq}\left(J_n\right) $ is generated by $  \mathcal{P}_{n-1} = \mathbb{M}_{n-1} $ as a vector space and $ \dim_\mathbb{C} \initial_{\preceq}\left(J_{n-1}\right) / \initial_{\preceq}\left(J_n\right) = n+1 $. However $ \sharp \mathcal{P}_{n-1} = n+2 $ by Lemma \ref{Pn-lemma} (1). Therefore there exists a non-trivial relation between monomials of $ \mathcal{P}_{n-1} $ in $ \initial_{\preceq}\left(J_{n-1}\right) / \initial_{\preceq}\left(J_n\right) $. Hence precisely one element $ \alpha \in \mathcal{P}_{n-1} $ is contained in $ \initial_{\preceq}\left(J_n\right) $; otherwise the existence of such relation contradicts to Lemma \ref{Macaulay'sTh}.
			
			We will see $ \alpha = p_{n-1} $. 
			
			Suppose that $ \alpha \in \mathcal{P}_{n-1} \setminus \left\lbrace \, p_{n-1} \, \right\rbrace $. Now there exists $ h \in J_n $ such that $ \lm{h} = \alpha $ since $ \alpha \in \initial_{\preceq}\left(J_n\right) $. Then $ h \in \ker \phi $ by Lemma \ref{phi-lemma} (3), and $ h =  (uv-1)f $ for some $ f \in S $ by Lemma \ref{phi-lemma} (1). Now $ f \in J_{n-1} $ by Lemma \ref{QuotDim-lemma} (2), and hence $ \lm{f} = \alpha - (1,1) \in \initial_{\preceq}\left(J_{n-1}\right) $. However this leads to a contradiction. Indeed, by the induction hypothesis, there exists $ \alpha' \in \mathcal{P}_{n-1} $ such that $ x^{\alpha'} $ divides $ x^{\alpha - (1,1)} $, that means $ \alpha \in \alpha' + (1,1) + \sigma_{\mathbb{Z}} \subset \alpha' + \sigma_{\mathbb{Z}} $. Now $ \alpha, \alpha' \in \mathcal{P}_{n-1} $ and this contradicts to Lemma \ref{Pn-lemma} (2).
			
			Thus $ \alpha = p_{n-1} $ and hence $ \alpha = p_n \in \initial_{\preceq}\left(J_n\right) $ by Lemma \ref{Pn-lemma} (2). Therefore $ \mathcal{P}_n \subset \initial_{\preceq}\left(J_n\right) $. This completes the proof.
		\end{proof}
	\end{Prop}

	\subsection{Non-regularity of $ C_{\mathbb{G}_n} $}
	
	Next our purpose is to show the non-regularity of $ C_{\mathbb{G}_n} $. We already saw that $ C_{\mathbb{G}_n} $ is a $ 2 $-dimensional cone (Lemma \ref{CG-lemma} (1)). Hence $ C_{\mathbb{G}_n} $ has two rays and we only have to determine them.
	
	Our strategy is as follows. First we choose a certain $ w \neq (0,0) $ from $ C_{\mathbb{G}_n} $. Next we find some $ \left(g,\alpha\right) \in \mathbb{G}_n $ such that $ \left(\alpha - \beta \right)\cdot w = 0 $ for some $ \beta \in \supp{g}\setminus \left\lbrace \, \alpha \, \right\rbrace $. Now $ \left(\alpha - \beta \right) \in C_{\mathbb{G}_n}^\vee $ according to the definition of $ C_{\mathbb{G}_n} $, and hence we can conclude that $ \mathbb{R}_{\geq 0}w $ is a ray of $ C_{\mathbb{G}_n} $.
	
	\begin{Lem}\label{supp(f)Dn-lemma}
		Let $ f \in J_n $ satisfy $ \text{lc}_{\preceq}\left( f \right) = 1 $. If $ \lm{f} \in \mathcal{P}_n $ and $ \supp{f} \setminus \left\lbrace \, \lm{f} \, \right\rbrace \subset \mathcal{D}_n $, then $ \left(f,\lm{f}\right) \in \mathbb{G}_n $. 
		
		\begin{proof}
			Let $ \alpha := \lm{f} $. Then $ \alpha \in \mathcal{P}_n = \mathbb{M}_n $ by Proposition \ref{lm-proposition} and hence there exists $ \left( g_\alpha, \alpha \right) \in \mathbb{G}_n $. Since $ \supp{f}\setminus\left\lbrace \, \alpha \, \right\rbrace \subset \mathcal{D}_n $, no monomial of $ \supp{f}\setminus\left\lbrace \, \alpha \, \right\rbrace $ is divisible by any monomial of $ \mathbb{M}_n $. This implies that $ \left\lbrace \, f \, \right\rbrace \cup \left\lbrace \, g\,\mid\, (g,\beta) \in \mathbb{G}_n \, \right\rbrace \setminus \left\lbrace \, g_\alpha \, \right\rbrace $ is a reduced Gr\"{o}bner basis of $ J_n $ with respect to $ \preceq $. By uniqueness of the reduced Gr\"{o}bner basis, we conclude $ f = g_\alpha $.
		\end{proof}
	\end{Lem}

	\begin{Prop}
		$ L_1 := \mathbb{R}_{\geq 0}(2,-1) $ is a ray of $ C_{\mathbb{G}_n} $.
		
		\begin{proof}
			Let $ w := (2,-1) $. Then $ C_{\mathbb{G}_n} $ contains $ w $ by Lemma \ref{w1-lemma}.
			
			As in Proposition \ref{n=1}, the reduced Gr\"{o}bner basis of $ J_1 $ contains
			$$
			g_1 := u^3v^4+u-4uv+2
			$$
			where $ \lm{g_1} = (3,4) = s_1 $. Now let $ g_n := \left( uv-1 \right)^{n-1}g_1 $.
			
			For any $ f,g \in S $, one can easily see that $ \initial_{w}\left(fg\right) = \initial_{w}\left(f\right)\initial_{w}\left(g\right) $. Now $ \initial_{w}\left( g_1 \right) = u^3v^4+u $ and $ \initial_{w}\left( \left( uv-1 \right)^{n-1} \right) = \left(uv\right)^{n-1} $. Therefore $ \initial_{w}\left( g_n \right) = \left(uv\right)^{n-1}\left(u^3v^4+u\right) $. Hence we have
			\begin{displaymath}
			\alpha_n := \left(uv\right)^{n-1}u^3v^4,\ \beta_n := \left(uv\right)^{n-1}u \in \supp{g_n}.
			\end{displaymath}
			
			Now $ \lm{g_n} = \alpha_n = (n-1)(1,1) + s_1  \in \mathcal{P}_n $ by Lemma \ref{Pn-lemma} (5). Moreover $ \supp{g_1}\setminus \left\lbrace \, \alpha_1 \, \right\rbrace \subset \mathcal{D}_1 $ and $ \alpha_1 \in \mathcal{D}_2 $. Hence, by Lemma \ref{Dn-lemma} (3) and $ \mathcal{D}_i \subset \mathcal{D}_{i+1} $, we have $ \supp{g_n}\setminus \left\lbrace \, \alpha_n \, \right\rbrace \subset \mathcal{D}_n $. Therefore $ \left( g_n, \alpha_n \right) \in \mathbb{G}_n $ by Lemma \ref{supp(f)Dn-lemma}.
			
			Then the following vector
			\begin{displaymath}
			\alpha_n - \beta_n = (3,4) - (1,0) = (2,4)
			\end{displaymath}
			satisfies $ (2,4)\cdot w = 0 $. Hence as we explained our strategy, $ L_1 $ is a ray of $ C_{\mathbb{G}_n} $.
		\end{proof}
	\end{Prop}

	\begin{Lem}\label{GiveAnotherRay-lemma}
		Let $ n \geq 2 $ be an integer. Then $ \mathbb{G}_n $ contains $ \left(g,\alpha \right) $ with the following conditions; if $ n $ is odd,
		\begin{displaymath}
		\alpha = q_n^{\frac{n-1}{2}}\ \text{and}\  r_{n-1}^{\frac{n-1}{2}} \in \supp{g};
		\end{displaymath}
		if $ n $ is even,
		\begin{displaymath}
		\alpha = p_n\ \text{and}\  s_{n-1} \in \supp{g}.
		\end{displaymath}
		
		\begin{proof}			
			Let $ n $ be even. Let $ \phi:\ S \rightarrow \mathbb{C}[\lambda^\pm] $ be the morphism in Lemma \ref{phi-lemma}, and 
			\begin{displaymath}
			\overline{f} := -\left( \lambda^{-1}-1 \right)^{\frac{n}{2}+1}\left( \lambda - 1 \right)^\frac{n}{2} \in \phi\left(J_{n}\right).
			\end{displaymath}
			One can easily check that $ \lambda^{-(\frac{n}{2}+1)},\ \lambda^{\frac{n}{2}} \in \supp{\overline{f}} $ and
			\begin{displaymath}
			\supp{\overline{f}} \subset C := \left\lbrace \, \lambda^{-(\frac{n}{2}+1)}, \lambda^{-\frac{n}{2}},\ldots,\lambda^{-1},1,\lambda,\ldots,\lambda^{\frac{n}{2}} \, \right\rbrace.
			\end{displaymath}
			Now Lemma \ref{Dn-lemma} (4) shows that any monomial in $ C \setminus\left\lbrace \, \lambda^{-(\frac{n}{2}+1)} \, \right\rbrace $ has a preimage by $ \phi $ in $ \mathcal{D}_n $, and in particular, $ \lambda^{\frac{n}{2}} $ has the only preimage $ s_{n-1} \in \mathcal{D}_n $. In addition, $ \lambda^{-(\frac{n}{2}+1)} $ has the preimage $ p_n $. Hence one can find a preimage $ f $ of $ \overline{f} $ such that
			\begin{displaymath}
			 p_n, s_{n-1} \in \supp{f}\ \text{with}\ \supp{f} \setminus \left\lbrace \, p_n \, \right\rbrace \subset \mathcal{D}_n.
			\end{displaymath}
		
			Since $ \phi(f) = \overline{f} \in \phi(J_{n}) $, there exists $ \Delta \in \ker(\phi) $ such that $ f + \Delta \in J_{n} $. Let $ \left\lbrace \, g_1,\ldots,g_t \, \right\rbrace $ be the reduced Gr\"{o}bner basis of $ J_{n} $ with respect to $ \preceq $. Then, by the division algorithm [D2; Appendix A, Theorem A.1.4], $ \Delta $ has the following expression;
			\begin{displaymath}
			\Delta = \sum_{i = 1}^t q_ig_i + r\ \text{where}\ \supp{r} \subset \mathcal{D}_{n}.
			\end{displaymath}
			Now $ r $ also satisfies $ g:= f + r \in J_{n} $ because $ g = \left(f + \Delta\right) - \sum q_ig_i $ is a difference of elements of $ J_n $. We will see that this $ g $ is the one we expects.
			
			It is clear that $ \supp{g} \setminus \left\lbrace \, p_n \,\right\rbrace \subset \mathcal{D}_n $ by the similar conditions of supports of $ f $ and $ r $. Moreover $ p_n \in \supp{g} $ since $ p_n \in \supp{f} \setminus \supp{r} $. Thus $ \lm{g} = p_n $ by Lemma \ref{Dn-lemma} (5) and hence $ \left(g, p_n\right) \in \mathbb{G}_n $ by Lemma \ref{supp(f)Dn-lemma}.
			
			Too see $ s_{n-1} \in \supp{g} $, it is enough to check $ s_{n-1} \notin \supp{r} $ since $ s_{n-1} \in \supp{f} $. The coefficient of $ s_{n-1} $ in $ r $ coincides with the coefficient of $ \phi\left(s_{n-1}\right) = \lambda^{\frac{n}{2}} $ in $ \phi(r) $ because of $ \supp{r} \subset \mathcal{D}_n $ and Lemma \ref{Dn-lemma} (4). However one can see $ \phi(r) = 0 $ as follows. Since $ \Delta \in \ker \phi $ and $ g_i \in J_n $, we have $ \phi(r) = \phi\left( \Delta - \sum q_ig_i \right) = 0 $ in $ \mathbb{C}[\lambda^\pm] / \phi\left(J_n\right) $. Therefore $ \phi(r) $ gives a linear relation between monomials of $ \phi\left(\supp{r}\right)\subset \phi\left( \mathcal{D}_n \right) $ in $ \mathbb{C}[\lambda^\pm] / \phi\left(J_n\right) $. Now $ \phi\left( \mathcal{D}_n \right) $ is described in Lemma \ref{Dn-lemma} (4), and Lemma \ref{phi-lemma} (2) shows that the linear relation must be trivial since $ \phi\left(\mathcal{D}_n\right) \subset B $. This implies $ \phi(r) = 0 $. Hence the coefficient of $ s_{n-1} $ in $ r $ is zero. Therefore $ s_{n-1} \notin \supp{r} $. Hence $ s_{n-1} \in \supp{g} $ and the case of even $ n $ is done.
			
			Let $ n $ be odd. By the case of even $ n $, there exists $ \left(h,p_{n-1}\right) \in \mathbb{G}_{n-1} $ such that $ s_{n-2} \in \supp{h} $.
			
			Let $ g := (uv-1)h \in J_n $. It is clear that $ \lm{g} = p_{n-1} + (1,1) = q_n^{\frac{n-1}{2}} \in \mathcal{P}_n $ by Lemma \ref{Pn-lemma} (4). One can easily see that $ \supp{g}\setminus\left\lbrace \, \lm{g} \, \right\rbrace \subset \mathcal{D}_n $ because of Lemma \ref{Dn-lemma} (3) with $ \supp{h}\setminus\left\lbrace \, p_{n-1} \, \right\rbrace \subset \mathcal{D}_{n-1} $ and $ p_{n-1} \in \mathcal{D}_n $. Thus Lemma \ref{supp(f)Dn-lemma} shows $ \left(g,q_n^{\frac{n-1}{2}}\right) \in \mathbb{G}_n $.
			
			Moreover $ r_{n-1}^{\frac{n-1}{2}} \in \supp{g} $. Indeed $ s_{n-2} \in \supp{h} $ and $ r_{n-1}^{\frac{n-1}{2}} = s_{n-2} + (1,1) $ by Lemma \ref{Pn-lemma} (4). Now $ \supp{h} \setminus \left\lbrace \, p_{n-1} \, \right\rbrace \subset \mathcal{D}_{n-1} $ and hence $ r_{n-1}^{\frac{n-1}{2}} \notin \supp{h} $. Therefore there is no cancellation at $ r_{n-1}^{\frac{n-1}{2}} $ between $ (uv-1)h = g $, and hence $ r_{n-1}^{\frac{n-1}{2}} \in \supp{g} $.
		\end{proof}
	\end{Lem}

	\begin{Prop}
		Let $ L_2 :=
		\begin{cases}
		\mathbb{R}_{\geq 0}(2n-2,-n+2) & \text{if}\ n\ \text{is odd,}\\
		\mathbb{R}_{\geq 0}(2n,-n+1) & \text{if}\ n\ \text{is even}.\\
		\end{cases}
		$
		
		Then $ L_2 $ is a ray of $ C_{\mathbb{G}_n} $.
		
		\begin{proof}
			Let $ l_n $ be the one in Lemma \ref{Dn-lemma} (6);
			$$ l_n :=
			\begin{cases}
			(2n-2,-n+2) & \text{if}\ n\ \text{is odd,}\\
			(2n,-n+1) & \text{if}\ n\ \text{is even.}
			\end{cases}
			$$
			Then $ L_2 = \mathbb{R}_{\geq 0}l_n $.
			
			First, let us show $ l_n \in C_{\mathbb{G}_n} $. By the definition of $ C_{\mathbb{G}_n} $, it is enough to check $ l_n \cdot (\alpha - \beta) \geq 0 $ for any $ \left(g,\alpha \right) \in \mathbb{G}_n $ and $ \beta \in \supp{g}\setminus \left\lbrace \, \alpha \, \right\rbrace $. Proposition \ref{lm-proposition} shows that $ \alpha \in \mathcal{P}_n $ and $ \beta \in \mathcal{D}_n $. In Lemma \ref{Dn-lemma} (6), we already saw that $ \Psi_n:\ \sigma_\mathbb{Z} \in a \mapsto l_n \cdot a \in \mathbb{R}_{\geq 0} $ satisfies
			\begin{displaymath}
			\max \Psi_n\left(\mathcal{D}_n\right) = \min \Psi_n\left(\mathcal{P}_n\right).
			\end{displaymath}
			Therefore $ l_n \cdot (\alpha - \beta) \geq 0 $ and hence $ l_n \in C_{\mathbb{G}_n} $ .
			
			On the other hand, by Lemma \ref{GiveAnotherRay-lemma}, $ \mathbb{G}_n $ contains $ \left(g,\alpha \right) $ such that
			\begin{align*}
			&\alpha = q_n^{\frac{n-1}{2}}\ \text{and}\  r_{n-1}^{\frac{n-1}{2}} \in \supp{g}\ \text{for odd}\ n,\\
			&\alpha = p_n\ \text{and}\  s_{n-1} \in \supp{g}\ \text{for even}\ n.
			\end{align*}
			Now let $ v_n := \begin{cases}
			q_n^{\frac{n-1}{2}} - r_{n-1}^{\frac{n-1}{2}} & \text{if}\ n\ \text{is odd,}\\
			p_n - s_{n-1} & \text{if}\ n\ \text{is even}
			\end{cases} $ then $ l_n \cdot v_n = 0 $ by Lemma \ref{Dn-lemma} (6). Therefore $ L_2 $ is a ray of $ C_{\mathbb{G}_n} $.
		\end{proof}
	\end{Prop}

	Now $ C_{\mathbb{G}_n} $ is completely described; $ C_{\mathbb{G}_n} $ is the $ 2 $-dimensional cone whose rays $ L_1 = \mathbb{R}_{\geq 0}(2,-1) $ and $ L_2 = \mathbb{R}_{\geq 0}l_n $ where
	$$ l_n :=
	\begin{cases}
	(2n-2,-n+2) & \text{if}\ n\ \text{is odd,}\\
	(2n,-n+1) & \text{if}\ n\ \text{is even.}
	\end{cases}
	$$
	One can easily check that $ l_n $ is the primitive ray generator of $ L_2 $.
	
	The following theorem is our main result.

	\begin{Th}\label{Non-regularCG-theorem}
		For any $ n > 0 $, $ \overline{\text{Nash}_n(X)} $ has a singular point of type $ A_1 $, and $ \text{Nash}_n(X) $ is singular.
	
		\begin{proof}
			To see the non-regularity of $ C_{\mathbb{G}_n} $, let $ N $ be the sublattice of $ \mathbb{Z}^2 $ generated by $ w := (2,-1) $ and $ l_n $. Then $ N \neq \mathbb{Z}^2 $ since
			\begin{displaymath}
			\det \left(\begin{matrix}
			w\\
			l_n
			\end{matrix}\right)
			= 2.
			\end{displaymath}
			Hence $ C_{\mathbb{G}_n} $ is non-regular. Moreover this calculation shows that the affine toric variety associated to $ C_{\mathbb{G}_n} $ is the $ A_1 $-singularity $ \left(z^2 -xy = 0 \right) \subset \mathbb{A}^3 $. Thus $ \overline{\text{Nash}_n(X)} $ has a singular point of type $ A_1 $ by Theorem \ref{GF-theorem}. Hence $ \text{Nash}_n(X) $ is also singular; otherwise $ \overline{\text{Nash}_n(X)} = \text{Nash}_n(X) $ and one has a contradiction.
		\end{proof}
	\end{Th}


\begin{thebibliography}{20}
		\bibitem[CLO]{CLO} Cox, D.\ A., Little, J., O'Shea, D.\ (2005).\ \textit{Using Algebraic Geometry}, 2nd ed., Graduate Text in Math.\ \textbf{185}, Springer.
		\bibitem[D1]{Duarte1} Duarte, D.\ \textit{Higher Nash blowup on normal toric varieties}, Journal of Algebra \textbf{418}: 110-128.
		\bibitem[D2]{Duarte2} Duarte, D.\ \textit{Nash modification on toric surfaces and higher Nash blowup on normal toric varieties}, Ph.\ D.\ thesis at University of Toulouse.
		\bibitem[H]{Hironaka} Hironaka, H.\ \textit{Resolution of singularities of an algebraic variety over a field of charactaristic zero I, II}, Ann.\ of Math.\ (2) \textbf{79}.
		\bibitem[N]{Nobile} Nobile, A. \textit{Some properties of the Nash blowing-up}, Pacific Journal of Math.\ \textbf{60}: 297-305.
		\bibitem[OZ]{OZ} Oneto, A., Zatini, E.\ \textit{Remarks on Nash blowing-up}, Rend.\ Sem.\ Mat.\ Univ.\ Torino \textbf{49}: 71-82.
		\bibitem[R]{Rob} Robbiano, L.\ \textit{Term Ordering on the Polynomial Ring}, Lecture Notes in Computer Science \textbf{204}, Springer: 513-517.
		\bibitem[S]{Sturmfels} Sturmfels, B.\ (1996).\ \textit{Gr\"{o}bner Bases and Convex Polytopes}, University Lecture Series \textbf{8}, Amer.\ Math.\ Soc.
		\bibitem[Y1]{Yasuda1} Yasuda, T.\ {\it Higher\ Nash\ blowups},\ Compositio\ Math.\ \textbf{143}: 1493-1510.
		\bibitem[Y2]{Yasuda2} Yasuda, T.\ {\it Universal\ flattening\ of\ Frobenius}, American Journal of Math.\ \textbf{134}, No.\ 2: 349-378.
	\end{thebibliography}
\end{document}